\documentclass[12pt]{article}
\usepackage{amsmath,amssymb}
\usepackage{pstricks}
\usepackage{pst-node}
\usepackage{latexsym}
\usepackage{color,graphicx}
\usepackage{wrapfig}
\newtheorem{Lemma}{Lemma}
\newtheorem{Proposition}[Lemma]{Proposition}

\newtheorem{Theorem}[Lemma]{Theorem}

\newtheorem{Corollary}[Lemma]{Corollary}
\newcommand{\con}{\rightarrow}

\newcommand{\cd}{\ \stackrel{d}{\rightarrow} \ }
\newcommand{\cp}{\ \stackrel{p}{\rightarrow} \ }

\newcommand{\GG}{\mbox{${\cal G}$}}
\newcommand{\FF}{\mbox{${\cal F}$}}

\newcommand{\eps}{\varepsilon}

\newcommand{\bX}{{\bf X}}

\newcommand{\bck}{\!\!\!}

\newcommand{\sfrac}[2]{{\textstyle\frac{#1}{#2}}}

\newcommand{\bY}{{\bf Y}}
\newcommand{\bx}{{\bf x}}
\newcommand{\by}{{\bf y}}
\newcommand{\bc}{{\bf c}}

\newcommand{\MC}{mul\-ti\-pli\-ca\-tive co\-a\-le\-scent}

\newcommand{\cvd}{l^2_{\mbox{{\footnotesize $\searrow$}}}}
\newcommand{\cvt}{l^3_{\mbox{{\footnotesize $\searrow$}}}}
\newcommand{\cvtwd}{\cvt \backslash \,\, \cvd}

\newcommand{\lrarrow}{\leftrightarrow}

\def\endpf{$\Box$}
\def\endrem{$\diamond$}
\newcommand{\Zb}{{\bar Z}}
\newcommand{\Yb}{{\bar Y}}
\newcommand{\Rb}{{\bar R}}

\newcommand{\len}{\rm len}

\begin{document}
\title{Eternal multiplicative coalescent is encoded by its L{\'e}vy-type processes}
\date{}
\author{Vlada Limic\thanks{
        UMR 8628, Département de Math\'ematiques,
        CNRS et Universit\'e Paris-Sud XI,
        91405 Orsay, France}\\
        e-mail: vlada.limic@math.u-psud.fr}

\maketitle
\vspace{-1cm}
\begin{abstract}
The multiplicative coalescent is a Markov process taking values in ordered $l^2$.  
It is a mean-field process in which any pair of blocks coalesces at rate proportional to
the product of their masses.
In Aldous and Limic (1998)  each extreme eternal version
$(\bX(t),- \infty < t < \infty)$ of the
multiplicative coalescent was described in three different ways.
One of these specifications matches
the (marginal) law of $\bX(t)$  to that of the ordered excursion lengths above past minima of $\{L_{\bX}(s) +ts, \,s \geq 0\}$, where
$L_{\bX}$ is a certain L{\'e}vy-type process which (modulo shift and scaling) has infinitesimal drift  $-s$ at time $s$.

Using a modification of the breadth-first-walk construction from Aldous (1997) and Aldous and Limic (1998), and some
 new insight from the thesis by Uribe (2007), this work settles an open problem (3) from Aldous (1997), in the more general context of Aldous and Limic (1998).
Informally speaking, $\bX$ is entirely encoded by $L_\bX$, and contrary to Aldous' original intuition, the evolution of time for $\bX$ does correspond to the linear increase in the constant part of the drift of $L_\bX$.
In the ``standard multiplicative coalescent'' context of Aldous (1997), this result was first
announced by Armend\'ariz in 2001, while its first published proof is due to Broutin and Marckert (2016), who simultaneously account  for the process of excess (or surplus) edge counts. 

The novel argument presented here is based on a sequence of relatively elementary observations.
Some of its components (for example, the new dynamic random graph construction via ``simultaneous'' breadth-first walks) are of independent interest, and 
may turn out useful for obtaining more sophisticated asymptotic results on near critical random graphs and related processes. 
\end{abstract}



\newpage
\section{Introduction}
\subsection{The \MC\ in 1997}
\label{sec1.1}
Recall that  the {\em multiplicative coalescent} takes values in the space of
collections
of blocks, where each block has mass in $(0,\infty)$, and informally evolves according to the following dynamics:
\begin{equation}
\label{merge}
\begin{array}{c}
\mbox{ each pair of blocks of mass $x$ and $y$ merges at rate $xy$ }\\
\mbox{into a single block of mass $x+y$.}
\end{array}
\end{equation}
It is well-known that the multiplicative coalescent is strongly connected to the (Erd\"os-R\'enyi \cite{erdren}) random graph, viewed in continuous time. 
More precisely, if $x_i$ is a positive integer for each $i$, we can represent each initial block as a collection of $x_i$ different particles of mass $1$, which have already merged in some specified arbitrary way. 
For each pair of particles $\{k,l\}$ let $\xi_{k,l}$ be an exponential
(rate $1$) random variable, and let $\xi_{k,l}$s be independent over $k,l$.
Let the graph evolve in continuous time according to the following mechanism: at time $\xi_{k,l}$ the particles $k$ and $l$ get connected by a (new) edge. Two different connected components merge at the minimal connection time of a pair of particles $(k,l)$, where $k$ is from one, and $l$ from the other component.
The mass of any connected component equals its number of particles.
Then the process of component masses evolves according to (\ref{merge}).

For a given initial state with a finite number of blocks (but block masses not necessarily integer-valued), it is easy to formalize
(\ref{merge}), e.g.~via a similar ``graph-construction'', in order to define a continuous-time finite-state Markov process.
Furthermore, Aldous \cite{aldRGMC} shows that one can extend the state space to include $l_2$ configurations.
More precisely, if $(\cvd,d)$ is  the metric space of infinite sequences
${\bf x} = (x_1,x_2,\ldots)$
with $x_1 \geq x_2 \geq \ldots \geq 0$
and $\sum_i x_i^2 < \infty$,
where 
$d({\bf x},{\bf y}) = \sqrt{\sum_i (x_i-y_i)^2}$, then the {\em multiplicative coalescent} is a Feller process on $\cvd$ (see \cite{aldRGMC} Proposition 5, or Section 2.1 in \cite{VLthesis} for an alternative argument), evolving according to description (\ref{merge}).
The focus in \cite{aldRGMC} was on 
the existence and properties of the \MC, as well as on the construction of a particular eternal version  
$(\bX^*(t), -\infty < t <\infty)$, called the {\em standard} \MC.
The standard version arises as a limit of
the classical random graph process 
near the phase transition (each particle has initial mass $n^{-2/3}$ and the random graph is viewed at times $n^{1/3}+O(1)$).
In particular, the marginal distribution $\bX^*(t)$ of $\bX^*$ was described in \cite{aldRGMC} as follows:
if $(W(s),0 \leq s < \infty)$ is standard Brownian motion
and 
\begin{equation}
 W^t(s) = W(s) - \sfrac{1}{2}s^2  + ts, \ s \geq 0, \label{defWt}
\end{equation}
%
and $B^t$ is its ``reflection above past minima"
\begin{equation}
 B^t(s) = W^t(s) - \min_{0 \leq s^\prime \leq s} W^t(s^\prime), \ s \geq 0, \label{defBt}
\end{equation}
then (see \cite{aldRGMC} Corollary 2) the ordered sequence of excursion (away from $0$) lengths of $B^t$ 
has the same distribution as $\bX^*(t)$.
Note in particular that the total mass
$\sum_i X^*_i(t)$ is infinite.

The author's thesis \cite{VLthesis} was based on a related question: are there any other eternal versions of the \MC, and, provided that the answer is positive, what are they? 
Paper \cite{EBMC} completely described the {\em entrance boundary}  of the \MC\ (or equivalently, the set of all of its extreme eternal laws).
The {\em extreme eternal laws} or {\em versions} are conveniently characterized by the property that their corresponding tail $\sigma$-fields at time $- \infty$ are trivial. 
Any (other) eternal versions must be a mixture of extreme ones, see e.g.~\cite{dyn78} Section 10.
Note that the word ``version'' is not used here in a the classical (Markov process theory) sense. 

\subsection{Characterizations of eternal versions in 1998}
The notation to be introduced next is inherited from \cite{EBMC}.
We write
$\cvt$ for the space of infinite sequences
${\bf c} = (c_1,c_2,\ldots)$
with $c_1 \geq c_2 \geq \ldots \geq 0$
and $\sum_i c_i^3 < \infty$.
For $\bc \in \cvt$,
let $(\xi_j, j \geq 1)$ be independent with exponential (rate $c_j$)
distributions and consider
\begin{equation}
 V^\bc (s) = \sum_j \left(c_j 1_{(\xi_j \leq s)} - c_j^2s \right)
, \ s \geq  0 . \label{defVc}
\end{equation}
We may regard 
$V^\bc$ as a L{\'e}vy-type process, where for each $x$ only the first jump of size $x$ is kept  (cf.~Section 2.5 of \cite{EBMC}, and Bertoin \cite{bertoin-levy} for background on L{\'e}vy processes).
It is easy to see
that $\sum_i c_i^3 < \infty$ is precisely  the condition for (\ref{defVc}) to yield a well-defined 
process (see also Section 2.1 of \cite{EBMC}).

Define the parameter space
$$
{\cal I} :=
\left( (0,\infty) \times (-\infty,\infty) \times \cvt\right) \cup \left( \{0\} 
\times (-\infty,\infty) \times {\cvtwd} \right).
$$
Now modify (\ref{defWt},\ref{defBt}) by defining,
for each $(\kappa, \tau, \bc) \in {\cal I}$,
\begin{equation}
\label{deftilW}
\widetilde{W}^{\kappa,\tau}(s) = \kappa^{1/2}W(s) +
\tau s - \sfrac{1}{2}\kappa s^2 , \ s \geq 0 
\end{equation}
\begin{equation}
W^{\kappa, \tau, \bc}(s) = \widetilde{W}^{\kappa, \tau}(s) + V^\bc (s), \ s \geq 0 \label{defWtc}
\end{equation}
\begin{equation}
 B^{\kappa, \tau,\bc}(s) = W^{\kappa, \tau,\bc}(s) - \min_{0 \leq s^\prime \leq s} W^{\kappa, \tau,\bc}(s^\prime), \ s \geq 0. \label{defBtc}
\end{equation}
So 
$B^{\kappa, \tau,\bc}(s)$
is again the reflected process with some set of (necessarily all finite, see Theorem \ref{TmainEBMC} below) excursions away from $0$.

Denote by $\hat{\mu}(y)$ the distribution of the constant process
\begin{equation}
\bX(t) = (y,0,0,0,\ldots), \ -\infty < t < \infty
\label{constant}
\end{equation}
where $y \geq 0$ is arbitrary but fixed. 

Let $\bX(t) = (X_1(t),X_2(t),\ldots)\in \cvd$ be the state of a particular eternal version of \MC. 
Then $X_j(t)$ is the mass of its $j$'th largest block at time $t$.
Write
\[ S(t) = S_2(t)= \sum_i X_i^2(t) , \mbox { and } S_3(t) = \sum_i X^3_i(t) .   \]

The main results of \cite{EBMC} are stated next. 
\begin{Theorem}[\cite{EBMC}, Theorems 2--4]
\label{TmainEBMC}
(a) For each
$(\kappa,\tau,\bc) \in {\cal I}$
there exists an eternal
\MC\ $\bX$ such that for each $-\infty < t < \infty$,
$\bX(t)$ is distributed as the 
ordered sequence of excursion lengths of $B^{\kappa, t - \tau, \bc}$. \\
(b)  Denote by $\mu(\kappa, \tau, \bc)$ the distribution of  $\bX$ from (a).
The set of extreme eternal \MC\ distributions is precisely\\
$\{ \mu(\kappa, \tau, \bc): \ (\kappa, \tau, \bc) \in {\cal I}\}$
$ \cup\ \{\hat{\mu}(y): 0 \leq y < \infty\}$.\\
(c) Let $(\kappa,\tau,\bc) \in {\cal I}$.
An (extreme) eternal \MC\ $\bX$ has distribution
$\mu(\kappa, \tau, \bc)$
if and only if
\begin{eqnarray}
|t|^3 S_3(t) & \con \kappa + \sum_j c_j^3 & \mbox{ a.s. as } t \con - \infty \label{Thyp1} \\
t + \frac{1}{S(t)} & \con \tau & \mbox{ a.s. as } t \con - \infty \label{Thyp2} \\
|t| X_j(t) & \con c_j & \mbox{ a.s. as } t \con - \infty, \ \forall j \geq 1. \label{Thyp3} 
\end{eqnarray}
\end{Theorem}

In terms of the above defined parametrization,  the Aldous \cite{aldRGMC} standard (eternal) \MC\ has
distribution $\mu(1,0, {\bf 0})$.
The parameters $\tau$ and $\kappa$ correspond to  time-centering and time/mass scaling respectively:
if $\bX$ has distribution
$\mu(1,0,\bc)$, then 
$\widetilde{\bX}(t) = \kappa^{-1/3} \bX(\kappa^{-2/3}(t - \tau))$
has distribution
$\mu(\kappa, \tau, \kappa^{1/3} \bc)$.
Due to (\ref{Thyp3}),  the components of $\bc$ may be interpreted as the relative
sizes of distinguished large blocks in the $t \to - \infty$ limit.

\subsection{The main results}
The rest of this work will mostly ignore the constant eternal \MC s.
For a given $(\kappa, \tau, \bc)\in {\cal I}$ we can clearly write 
$
W^{\kappa, t-\tau, \bc}(s) = W^{\kappa, -\tau, \bc}(s) +t s , \,s \geq 0.
$
The  L{\'e}vy-type process $W^{\kappa, -\tau, \bc}$ is particularly important for this work. As we will soon see,  $W^{\kappa, -\tau, \bc}$ matches $L_{\bX}$ from the abstract, as soon as $\bX$ has law $\mu(\kappa, \tau, \bc)$.
 
As noted in \cite{EBMC} and  in \cite{aldRGMC} beforehand, at the time there was no appealing intuitive explanation of why
excursions of a stochastic process would be relevant in describing the marginal laws in Theorem \ref{TmainEBMC}(a).
One purpose of this work is to offer a convincing explanation (see Proposition \ref{Pcoro} below and in Section \ref{S:BFWUD}, Lemma \ref{Lsimultcon} in Section \ref{S:SL} and Lemma \ref{LjointconX} in Section \ref{S:C}).
Furthermore, open problem (3) of \cite{aldRGMC} asks about the existence of a two parameter (non-negative) process $(B^t(s),\,s \geq 0, \,  t \in \mathbb{R})$ such that the excursion (away from $0$)  lengths of $(B^{\,\cdot}(s),\,s\geq 0)$ evolve as $X^*(\cdot)$.
The statement of this problem continues by offering an intuitive explanation for why 
 $\{{\rm reflected}(W^{1, 0, {\bf 0}}(s) + ts),\ s\geq 0,\ t\geq 0\}$ should not
be the answer to this problem.  
Aldous' argument is more than superficially convincing, but the striking reality is that the simplest of guesses is actually not too naive to be true.
Armend\'ariz \cite{armen_thesis} (as well as \cite{armen05}) obtained but never published this result,
and Broutin and Marckert \cite{bromar15} recently derived it via a different technique, while considering in addition the excess-edge data in agreement with \cite{aldRGMC} (thus improving on the Armend\'ariz claim).

This is not the only surprise.
Popular belief judges the breadth-first-walk construction, on which \cite{aldRGMC,EBMC} reside, as  ``inadequate''  and the main reason for the just described ``confusion'' in the statement of \cite{aldRGMC}, open problem (3). One of the main points of this work is to show the contrary. Indeed, a modification of the original breadth-first-walk from \cite{aldRGMC,EBMC}, combined with 
a rigorous formulation (see Proposition \ref{PUriArm}) of Uribe's \cite{uribe_thesis} graphical interpretation of the  {Armend\'ariz' representation} \cite{armen_thesis,armen05}, one arrives to the following claim of independent interest, here stated for readers' benefit in the simplest (purely) homogenous setting.\\
{\bf Claim (Proposition \ref{Pcoro}, special case)}
{\em Suppose that $x_1=x_2=\ldots=x_n=1$, for some $n\in \mathbb{N}$, and define for $q>0$
$$
Z^q(s):= \sum_{i=1}^n  1_{(\xi_i\, \leq \ qs)} -s , \ s\geq 0,
$$
where $\xi_i$, $i=1.\ldots,n$ is a family of i.i.d.~exponential (rate $1$) random variables. 
For each $q>0$, let ``blocks'' be the finite collection of excursions (above past minima) of $Z^q$, and for each block let its mass be the corresponding excursion length.
Set $X(0)$ be the configuration of $n$ blocks of mass $1$, and for $q>0$ let $X(q)$ be the configuration of blocks at time $q$ (if needed, list the blocks in the mass non-increasing order, and append infinitely many $0$s). Then $(X(q),\,q\geq 0)$ is a continuous-time Erd\"os-R\'enyi random graph, evolving according to (\ref{merge}).
}

To the best of our knowledge, even the ``static'' statement that matches the law of $X(q)$ to the law of the continuous-time homogeneous Erd\"os-R\'enyi random graph for each time $q$ separately, was not previously recorded (even though the analysis of \cite{uribe_thesis}, on pages 111-112 is implicitly equivalent). To have a glimpse at the power of this approach, the reader is invited to fix $c>1$, and consider the asymptotic behavior (as $n\to \infty$) of a related process $Z^{(1/n,\ldots,1/n), cn}(s):=\sum_{i=1}^n  \frac{1}{n}1_{(n\xi_i\, \leq \ cns )} -s$ , $s\geq 0$,
(see Section \ref{S:SBFW} or Proposition \ref{Pcoro} for notation)
in order to determine (in a few lines only) the asymptotic size of the giant component in the supercritical regime.
In addition, the simultaneous breadth-first walks framework allows for a particularly elegant treatment of surplus edges, postponed to \cite{edgsur}.

The analysis similar to that of \cite{EBMC} (to be done in Sections \ref{S:SL} and \ref{S:C}) now yields:
\begin{Theorem}
\label{Tmain}
Fix a L{\'e}vy-type process $W^{\kappa, -\tau, \bc}$, and for any $t\in (-\infty,\infty)$ define 
\[
W^{\kappa, t-\tau, \bc}(s):= W^{\kappa, -\tau, \bc}(s) +ts, \ s\geq 0.
\]  
Let $B^{\kappa, t-\tau, \bc}$ be defined as in (\ref{defBtc}).
For each $t$, let $\bX(t)=\bX^{\kappa,\tau,\bc}(t)$ be the infinite vector of ordered excursion lengths of $B^{\kappa, t-\tau, \bc}$ away from $0$. 
Then $(\bX(t), t\in (-\infty,\infty))$ is a c\`adl\`ag realization of $\mu(\kappa, \tau, \bc)$. 
\end{Theorem}

\subsection{Further comments on the literature and related work}
For almost two decades the only stochastic merging process widely studied by probabilists was the (Kingman) coalescent \cite{kin82,kin82a}. 
Starting with Aldous \cite{aldous_survey,aldRGMC}, and Pitman \cite{pit99},  Sagitov \cite{sag99}, and Donnelly and Kurtz \cite{dk99}, the main-stream probability research on coalescents was much diversified.

The Kingman coalescent  and, more generally, the mass-less (exchangeable) coalescents of \cite{pit99,sag99,dk99}
mostly appear in connection to the mathematical population genetics,  as universal (robust) scaling limits of genealogical trees 
(see for example \cite{moehle sagitov,schweinsberg_xi,blg1,schdur}, or a survey \cite{ensaios}).\\
The standard multiplicative coalescent is the universal scaling limit of numerous stochastic
 (typically combinatorial or graph-theoretic) homogeneous (or symmetric) merging-like models \cite{aldRGMC,addbrogol,aldpit,bhamidietal1,bhamidietal2,bhahof1,riordan,turova}. 
The ``non-standard'' eternal extreme laws from \cite{EBMC} are also scaling limits of inhomogeneous random graphs and related processes under appropriate assumptions \cite{EBMC,bhahof2,bhahof3}.

The two nice graphical constructions for coalescents with masses were discovered early on: by Aldous in \cite{aldRGMC} for the multiplicative case, and almost simultaneously by Aldous and Pitman \cite{dajp97sac} for the additive case 
(here any pair of blocks of mass $x$ and $y$ merges at rate $x+y$). 
The analogue of \cite{EBMC} in the additive coalescent case is again due to Aldous and Pitman \cite{dajp97ebac}.
No nice graphical construction for another (merging rate) coalescent with masses seems to have been found since. 
For studies of stochastic coalescents with general kernel see Evans and Pitman \cite{evapit} and Fournier \cite{four1,four2}.
Interest for probabilistic study of related Smoluchowski's equations (with general merging kernels) was also incited by \cite{aldous_survey}, see for example Norris \cite{norris}, Jeon \cite{jeon}, then Fournier and Lauren\c{c}ot \cite{fourlaur1,fourlaur2} and Bertoin \cite{bertoin_smolu}
for more recent, and Merle and Normand \cite{mernor1,mernor2} for even more recent developments. 
All of the above mentioned models are {\em mean-field}. 
See for example \cite{ls,barethveb,free15} for studies of (mass-less) coalescent models in the presence of spatial structure.

As already mentioned, Broutin and Marckert \cite{bromar15} obtain Theorem \ref{Tmain} in the standard \MC\ case, via {\em Prim's algorithm} construction invented for the purpose of their study, and notably different from the approach presented here.  
Before them Bhamidi et al.~\cite{bhamidietal1,bhamidietal2} proved f.d.d.~convergence 
for models similar to Erd\"os-R\'enyi random graph. 
For the standard additive coalescent, analogous results were obtained rather early by Bertoin \cite{bertoin_frag,bertoin_additive} and 
Chassaing and Louchard \cite{chalou}, and are rederived in \cite{bromar15}, again via an appropriate Prim's algorithm representation. 

In parallel to and independently from the research presented here, both Martin and R\'ath \cite{marrat_prep} and Uribe Bravo \cite{uribe_wip} have been studying closely related models and questions.
Their approaches seem to be quite different from the one taken here, with some notable similarities.
James Martin and Bal\'azs R\'ath \cite{marrat_prep} introduce a 
coalescence-fragmentation model called the {\em \MC\ with linear deletion (MCLD)}.
Here in addition to (and independently of) the multiplicative coalescence, each component is permanently removed from the system at a rate proportional to its mass (this proportionality parameter is denoted by $\lambda$).  
In the absence of deletion (i.e.~when $\lambda=0$), their ``tilt (and shift) operator'' representation of the MCLD leads to an alternative proof of Theorem \ref{Tmain}, sketched in detail in \cite{marrat_prep}, Section 6.1 (see \cite{marrat_prep}, Corollary 6.6). 
Further comments on links and similarities to \cite{marrat_prep} will be made along the way, most frequently in Section \ref{S:UD}.
Ger\'onimo Uribe \cite{uribe_wip} relies on a generalization of the construction from \cite{bromar15}, explains its links to Armend\'ariz' representation, and works towards another derivation of Theorem \ref{Tmain}.

The arguments presented in the sequel are elementary in part due to direct applications of a non-trivial result from \cite{EBMC},  Section 2.6 (depending on \cite{EBMC},  Section 2.5)
in Section \ref{S:C} (more precisely, Corollary \ref{C:immediate}). 
In comparison, (a) \cite{bromar15} also rely on the convergence results of \cite{aldRGMC} in the standard \MC\ setting, as well as additional estimates proved in \cite{addbrogol}, 
and (b) the analysis done in \cite{marrat_prep}, Sections 4 and 5 seems to be a formal analogue of that in \cite{EBMC}, Sections 2.5-2.6.
%

The present approach to Theorem \ref{Tmain} is of independent interest even in the standard \MC\ setting (where Section \ref{S:SL} would simplify further, since $\bc={\bf 0}$, and already Lemma 8 from \cite{aldRGMC} would be sufficient for making conclusions in Section \ref{S:C}).
In addition, it may prove useful for continued analysis of the \MC s, as well as various other processes in the \MC\ ``domain of attraction''.

The reader is referred to Bertoin \cite{bertoin-fragcoal} and Pitman \cite{pitman-stflour}
for further pointers to stochastic coalescence literature, and to Bollobas \cite{bol_book} and Durrett \cite{durrett_rgd} for the random graph theory and literature. 

\bigskip
The rest of the paper is organized as follows: Section \ref{S:SBFW} introduces the simultaneous breadth-first walks and explains how they are linked to the (marginal) law of the \MC, and the original breadth-first walks of \cite{aldRGMC,EBMC}.
Section \ref{S:UD} recalls Uribe's diagrams and includes Proposition \ref{PUriArm}, that connects the diagrams to the \MC.
In Section \ref{S:BFWUD} the simultaneous BFWs and Uribe's diagrams are linked, and as a result an important conclusion is made in Proposition \ref{Pcoro} (the generalized version of the claim which precedes Theorem \ref{Tmain}). 
All the processes considered in Sections \ref{S:SBFW}--\ref{S:BFWUD} have finite initial states.
Section \ref{S:SL} serves to pass to the limit where the initial configuration is in $\cvd$. 
The similarities to and differences from \cite{EBMC} are discussed along the way.
Theorem \ref{Tmain} is proved in Section \ref{S:C}. 
Several questions are included in Section \ref{S:Q} (the reader is also referred to the list of open problems given at the end of \cite{aldRGMC}).

{\bf Acknowledgements.}  
The author did not engage in thinking about \MC s for at least fifteen years, and it was undoubtedly the joint effort of Nicolas Fournier, Mathieu Merle, Justin Salez, Jean Bertoin and Ger\'onimo Uribe Bravo that brought her interest back.
Were it not for the spur of discussions with Nicolas, and Mathieu and Justin in the Spring of 2014, the conversation with Jean would have not taken place. Jean generously told the author of a somewhat unusual history of this problem, and encouraged her to contact Ger\'onimo. Ger\'onimo gracefully provided her with a copy of the relevant chapter of his thesis, along with several helpful remarks. 
The author is immensely grateful to all these colleagues for their show of trust and support at the moment where it was likely that her research could be stalled for an unpredictably long period of time.
She benefited in addition from past exchanges with Ger\'onimo Uribe, and with Nicolas Broutin and Jean-Fran\c{c}ois Marckert, that accelerated the writing of this article, and the recent exchanges with James Martin and Bal\'azs R\'ath who
provided several very useful suggestions for improved presentantion. 

\section{Simultaneous breadth-first walks}
\label{S:SBFW}
This section revisits the Aldous' breath-first walk construction of the \MC\ started from a finite vector $\bx$ from \cite{aldRGMC}, with two important differences (or modifications), which will be described along the way. 

Recall that ``breadth-first'' refers here to the order in which the vertices of a given connected graph (or one of its spanning trees) are explored. Such exploration process starts at the root, visits all of its children (these vertices become the 1st generation), then all the children of all the vertices from the 1st generation (these vertices become the 2nd generation), then all the children of  the 2nd generation,
and keeps going  until all the vertices (of all the generations) are visited, or until forever (if the tree is infinite).

Refer to $\bx=(x_1,x_2,x_3\ldots)\in \cvd$ as {\em finite}, if  for some $i\in \mathbb{N}$ we have $x_i=0$. Let the {\em length} of $\bx$ be the number $\len(\bx)$ of non-zero coordinates of $\bx$.
Fix a finite initial configuration $\bx\in \cvd$.
For each $i\leq \len(\bx)$ let $\xi_i$ have exponential (rate $x_i$)
distribution, independently over $i$. 

Given $\xi$, simultaneously for all $q>0$, 
we construct the (modified)
{\em breadth-first walk}
associated with $\bX(q)$ started from $\bX(0)=\bx$ at time $0$.
This simultaneity in $q$ is a {\em new feature} with respect to \cite{aldRGMC,EBMC}. 

The sequence $(\xi_i)_{i\leq \len(\bx)}$ will be used both for size-biased picking of the connected components, and for finding the merger events between the blocks. 
Fix $q>0$, and consider the sequence  $(\xi_i/q)_{i\leq \len(\bx)}$. Let us introduce the abbreviation $\xi_i^q:=\xi_i/q$.
The order statistics of $(\xi^q_i)_{i\leq \len(\bx)}$ are  $(\xi^q_{(i)})_{i\leq \len(\bx)}$. Define 
\begin{equation}
\label{defZbxq}
Z^{\bx,q}(s):= \sum_{i=1}^{\len(\bx)} x_i 1_{(\xi^q_i\, \leq \ s)} -s = \sum_{i=1}^{\len(\bx)} x_{(i)} 1_{(\xi^q_{(i)}\, \leq \ s)} -s, \ s\geq 0, \ q>0.
\end{equation}
In words, $Z^{\bx,q}$ has a unit negative drift and successive positive jumps, which occur precisely at times $(\xi^q_{(i)})_{i\leq \len(\bx)}$, and where the $i$th  successive jump is of magnitude $x_{(i)}$.

Here is the first important observation.
For each $q$, the \MC\ started from $\bx$ and evaluated at time $q$ can be constructed in parallel to $Z^{\bx,q}$ via a breadth-first walk coupling, similar to the one from \cite{aldRGMC,EBMC}. The interval $F_1^q:=[0,\xi^q_{(1)}]$ is  the first ``load-free'' period. 
Set $J_0:=\{1,2,\ldots,\len(\bx)\}$. 
At the time of the first jump of $Z^{\bx,q}$ we record
$$
\pi_1:= i \mbox{ if and only if } \xi_i=\xi_{(1)}, \mbox{ and }  J_1:=J_0\setminus \{\pi_1\},
$$
so that $\pi_1$ is the index of the first size-biased pick from $(\bx_i)_{i=1}^{\len(\bx)}$ using $\xi$s (or equally, $\xi^q$s). 
Furthermore, let us define for $l\leq \len(\bx)$
$$
\pi_l:=i \mbox{ if and only if } \xi_i=\xi_{(l)}, \ l\in \{1,\ldots, \len(\bx)\},  \mbox{ and $J_l:= J_{l-1}\setminus \{\pi_l\}$}.
$$ 
In this way, $(x_{\pi_1},x_{\pi_2},\ldots, x_{\pi_{\len(\bx)} } )$  is the size-biased random ordering of the initial non-trivial block masses.
As already noted,
the random permutation $\pi$ does not depend on $q$.

Let $\FF^q_s:=\sigma\{\{\{\xi^q_i>u\}:\, i\in J_0\},\, u\leq s\}$. Then $\FF^q=\{\FF^q_s,\,s\geq 0\}$ is the filtration generated by the arrivals of $\xi^q$s. 
Due to elementary properties of independent exponentials, it is clear that the above defined process
$Z^{\bx,q}$ is a continuous-time Markov chain with respect to $\FF^q$.
Indeed, given $\FF^q_{s}$, the (residual) clocks $\xi^q_i-s$ are again mutually independent, and moreover on the event $\{\xi^q_i>s\}$ we clearly have
$
P(\xi^q_i-s>u|\FF^q_s)  = e^{-x_i q u} = P(\xi_i^q>u) .
$
Furthermore, $\xi^q_{(1)}$ is a finite stopping time with respect to $\FF^q$ and  
\begin{equation}
\label{Elaweq}
P(\xi_i^q- \xi^q_{(1)}>u|\FF^q_{\xi_{(1)}}) 1_{(i\in J_1)} = e^{-x_i q u} 1_{(i\in J_1)}= P(\xi^q_i>u) 1_{(i\in J_1)}.
\end{equation}
Let $I_0=\emptyset$ and $I_1:= (\xi^q_{(1)}, \xi^q_{(1)} +x_{\pi_1}]$. Note that the length of the interval $I_1$ is the same (positive) quantity $x_{\pi_1}$ for all $q>0$.
During the time interval $I_1$ the dynamics ``listens for the children of $\pi_1$''. More precisely, if for some $j$ we have $\xi^q_j\in I_1$, or equivalently, if $\xi^q_j-  \xi^q_{(1)} \leq x_{\pi_1}$, we can interpret this as 
\[
\mbox{ edge } j \lrarrow \pi_1 \mbox{ appears before time } q \mbox{ in the \MC.}
\]
Indeed, as argued above, $P(\xi^q_j-  \xi^q_{(1)} > x_{\pi_1}|\FF^q_{\xi^q_{(1)}})=e^{-q x_j x_{\pi_1}}$, and this is precisely the multiplicative coalescent probability of the $j$th and the $\pi_1$st block not merging before time $q$.  

For any two reals $a<b$ and an interval $[c,d]$ where $0\leq c <d$, define the concatenation
$$
 (a,b]\oplus [c,d]:=(a+c,b+d].
$$
Recall that $I_1=(\xi^q_{(1)}, \xi^q_{(1)} +x_{\pi_1}]$, and define $N_1$ to be the number of $\xi^q$s that rung during $I_1$ (this is the size of the 1st generation in the exploration process).
For any $l\geq 2$ define recursively: if $I_{l-1}$ is defined
\begin{equation}
\label{ErecurIl}
I_l^q\equiv I_l:= \left\{\begin{array}{ll}
 I_{l-1} 
\oplus [0,x_{\pi_l}], & \mbox{provided }\frac{\xi_{(l)}}{q} \in   I_{l-1}\\
{\rm undefined}, & \mbox{otherwise}
\end{array}
 \right., 
\end{equation}
and if $I_l$ is defined in (\ref{ErecurIl}), let 
\begin{equation}
\label{ErecurNl}
N_l^q \equiv N_l:= \mbox{ the number of $\xi^q$s that rung during $I_l$},
\end{equation}
and otherwise let $N_l$ be (temporarily) undefined.
Since $\xi^q$s decrease in $q$,
the intervals $I^q_\cdot$ defined in this (coupling) construction do vary over $q$ (their endpoints decrease in $q$), but all of their lengths are constant in $q$. In fact, if defined, $I_l$ equals 
$(\xi_{(1)}^q, \xi_{(1)}^q+\sum_{m=1}^l x_{\pi_m}]$.
We henceforth abuse the notation and mostly omit the superscript $q$ when referring to $I$s or $N$s. 

During each $I_l\setminus I_{l-1}$ the coupling dynamics ``listens for the children of  $\pi_l$'', among all the $\xi^q$s  which have not been heard before (i.e.~they did not ring during $I_{l-1}$). 
If $I_l$ is defined in  (\ref{ErecurIl}), the set of children of $\pi_l$ in the above breadth-first order is precisely $J_{N_{l-1}} \setminus J_{N_l}$, which will be empty if and only if $N_l=N_{l-1}$.
The same memoryless property of exponential random variables as used above (e.g.~in (\ref{Elaweq})) ensures that  
\begin{eqnarray}
& & P(\xi^q_k \in I_l \setminus I_{l-1}|\FF_{\xi_{(1)}^q+x_{\pi_1}+\cdots + x_{\pi_{l-1}}}) 1_{(k \in J_{N_{l-1}})} = \label{Elaweq1}\\
& & P(k\in J_{N_{l-1}}\setminus J_{N_l} |\FF_{\xi_{(1)}^q+x_{\pi_1}+\cdots + x_{\pi_{l-1}}}) 1_{(k \in J_{N_{l-1}})} = e^{-qx_k x_{\pi_l}} 1_{(k \in J_{N_{l-1}})} \ \ \mbox{a.s.}\nonumber
\end{eqnarray}
Due to independence of $\xi$s, the residual clocks have again the (conditional) multi-dimensional product law.
So for each $l$, the set of children of  $\pi_l$ equals in law to the set of blocks which are connected by an edge to the $\pi_l$th block in the \MC\ at time $q$, given that they did not get connected by an edge (before time $q$) to any of the previously recorded blocks $\pi_1,\ldots,\pi_{l-1}$. 

The above procedure may (and typically will) stop at some $l_1\leq \len(\bx)$, due to $\xi^q_{(l_1)}$ not  arriving in  $I_{l_1-1}$. This will happen if and only if the whole connected component of the $\pi_1$st initial block  (in the \MC, evaluated at time $q$) was explored during $ I_{l_1-1}$, and the $\pi_{l_1-1}$st initial block was its last visited ``descendant'', while the rest of the graph was not yet ``seen'' during $s\in F_1^q \cup  I_{l_1-1}$. Indeed, if $a_1=\xi^q_{(1)}$ and $b_1=\xi^q_{(1)} + x_{\pi_1} + \ldots + x_{\pi_{l_1-1}}$, it is straight-forward to see that
\begin{equation}
\label{Eexcur}
 Z^{\bx,q}(s)> Z^{\bx,q}(a_1)=Z^{\bx,q}(b_1),  \ \forall s\in (a_1,b_1).
\end{equation}
In words, the interval $Cl(I_{l_1-1}) =[a_1,b_1]$ is an {\em excursion} of $Z^{q,\bx}$ {\em above past minima} of length  $b_1-a_1= x_{\pi_1} + \ldots + x_{\pi_{l_1-1}}$, which is the total mass of the first (explored) spanning tree in the breadth-first walk. 
Due to (\ref{Elaweq},\ref{Elaweq1}) and the related observations made above, this (random) tree matches the spanning tree of the connected component in the coupled \MC, evaluated at time $q$. This (first) spanning tree is rooted at $\pi_1$ (cf.~Figures 1 and 2) for all $q>0$. It will be clear from construction, that the roots of subsequently explored spanning trees can (and inevitably do) change at some $q>0$.

The next interval of time $F_2^q:= (\xi^q_{(1)} + x_{\pi_1}+\cdots + x_{\pi_{l_1-1}},\xi^q_{(l_1)}]$ is again ``load-free''  for the breadth-first walk.
Repeating the above exploration procedure starting from $\xi^q_{(l_1)}$  amounts to defining $I_{l_1}^q\equiv I_{l_1}:= (\xi^q_{(l_1)}, \xi^q_{(l_1)}+x_{\pi_{l_1}}]$ and listening for the children of $\pi_{l_1}$st block during $I_{l_1}$, and then running the recursion (\ref{ErecurIl},\ref{ErecurNl})  for $l\geq l_1+1$ until it stops, which occurs when all the vertices (blocks) of the second connected component are explored. 
This explorative coupling construction continues until all the initial blocks of positive mass are accounted for, or equivalently until $\xi^q_{(\len(\bx))}$.  
Clearly no $\xi$ can ring during $I_{\len(\bx)}\setminus I_{\len(\bx)-1}$ (which is open on the left), and $Z^{\bx,q}$ continues its evolution as a deterministic process (line of slope $-1$) starting from the left endpoint of $I_{\len(\bx)}$.

Figure 2 illustrates the just described coupling.
Each excursion of $Z^{\bx,q}$ above past minima corresponds uniquely to a connected component in the coupled \MC\ evaluated at time $q$.
It is clear from (\ref{ErecurIl},\ref{ErecurNl}) that the order of blocks visited within any given connected component is breadth-first. Note as well that the connected components are explored in the size-biased order. Indeed, the fact that the initial block of the next component to be explored is picked in a size-biased way, with respect to their individual masses, induces {\em size-biasing} of connected components (again with respect to mass) in the \MC\ at time $q$.

Let us fix some time $q$.
Figure 1 (without the vertical dashed lines and their labels) is a duplicate of \cite{EBMC}, Figure 1.
In the current notation  $\bx = (1.1,0.8,0.5,0.4,0.4,0.3,0.2,0,$
$0,\dots)$ so that $\len(\bx)=7$, and $x_{\pi_i}\equiv v(i)$. 
For the breath-first walks in \cite{aldRGMC,EBMC} the leading block in each component does not correspond to a jump of the walk, while
every non-leading block can be uniquely matched to a jump of the walk.
The time of this jump has exponential  (rate $q \cdot \text{ mass of the block}$) distribution,
and its size equals the mass of the non-leading block. All these exponential jumps are mutually independent.   
In particular, for the $i$th block, provided it is not leading, we can use $\xi_i^q=\xi_i/q$ in the construction of  Aldous' breadth-first walk from \cite{EBMC}. 

The here added vertical dashed lines and labels illustrate the link between the two breadth-first walk constructions (see also Figure 2 and the explanations provided below it). In particular, the first jump of Aldous' breadth-first walk happens at time $\xi_{(2)}^q - \xi_{(1)}^q$, the second one happens at time $\xi_{(3)}^q - \xi_{(1)}^q$ and this continues until the first component is exhausted. The next jump happens when the next non-leading block is encountered.

\vspace{0.2cm}
\includegraphics[width=350pt,height=340pt]{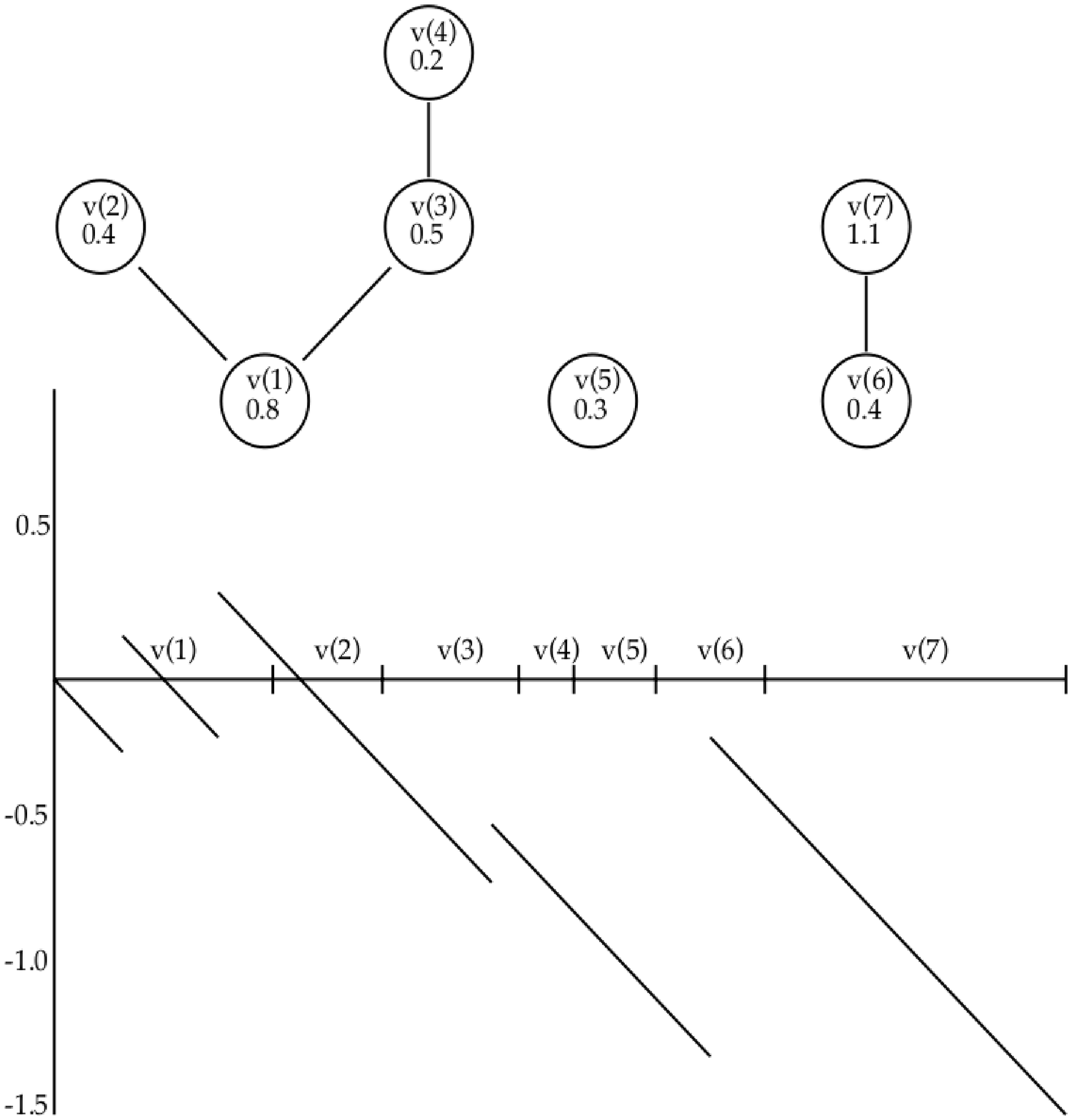}

   \psset{linewidth=0.4pt}
\psset{xunit=2.9cm,yunit=3.045cm}
   \pspicture(-0.353,-2.8)(10,-2.8)
\psline[linestyle=dashed](0.25,-2.7)(0.25,-0.5) 
\rput[t](0.25,-2.7){$\xi_{(2)}^q-\xi_{(1)}^q$}
\psline[linestyle=dashed](0.6,-2.4)(0.6,-0.5)  
\rput[t](0.6,-2.4){$\xi_{(3)}^q-\xi_{(1)}^q$}
\psline[linestyle=dashed](1.6,-2.5)(1.6,-0.5)  
\rput[t](1.6,-2.5){$\xi_{(4)}^q-\xi_{(1)}^q$}
\psline[linestyle=dashed](2.39,-2.7)(2.39,-0.5)  
\rput[t](2.4,-2.7){$\xi_{(7)}^q-\xi_{(6)}^q+ \sum_{i=1}^5 x_{\pi_i}$}
\endpspicture

\vspace{0.8cm}
\centerline{Figure 1}

\smallskip
From Figure 1 we can read off that $\pi$ belongs to $\{(2,4,3,7,6,5,1)$, $(2,5,3,7,6,4,1)\}$.
However, $\xi_{(1)}=\xi_2$, $\xi_{(5)}=\xi_6$ and $\xi_{(6)}=\xi_5$ (or $\xi_{(6)}=\xi_4$, depending on $\pi$) are not observed. 
In the simultaneous
 breadth-first walks construction the $\xi_{(1)}^q$ (here $\xi_2$), $\xi_{(5)}$ (here $\xi_6$) and $\xi_{(6)}$ also influence the walk. 

  \psset{linewidth=0.4pt}
\psset{xunit=2.9cm,yunit=3.045cm}
   \pspicture(0,-1)(4.4,1.2)
   \psset{linestyle=solid}
  \psline{-|}(0.1,0)\psline{-|}(0.1,0)(2.3,0)\psline{-|}(2.3,0)(2.8,0)\psline{->}(2.3,0)(4.4,0)
   \psline{->}(0,1)\psline{-}(0,-1.0)
   \uput[l](4.35,-0.07){\small{$s$}}
    \psset{linewidth=0.6pt}
    \psline{-o}(0,0)(0.1,-0.1)            
    \psline{*-o}(0.1,0.7)(0.35,0.45) 
    \psline{*-o}(0.35,0.85)(0.7,0.5)  
    \psline{*-o}(0.7,1)(1.7,0)   
    \psline{*-o}(1.7,0.2)(2.3,-0.4)   
    \psline{*-o}(2.3,-0.1)(2.8,-0.6)   
    \psline{*-o}(2.8,-0.2)(3,-0.4)        
\psline{*-}(3,0.7)(4.2,-0.5)         
    \psset{dotsize=0.5mm}
    \psdots(4.25,-0.55)(4.3,-0.6)(4.35,-0.65) 
\psline[linewidth=0.8pt,linecolor=gray]{-*}(0,-0.8)(0.1,-0.8)
\psline[linewidth=0.8pt,linecolor=gray]{-*}(2,-0.8)(2.3,-0.8)
\psline[linewidth=0.8pt,linecolor=gray]{-*}(2.6,-0.8)(2.8,-0.8)
\psline[linewidth=0.8pt,linecolor=blue]{(-}(0.1,-0.8)(0.9,-0.8)
\psline[linewidth=0.8pt,linecolor=blue]{(-}(0.9,-0.8)(1.3,-0.8)
\psline[linewidth=0.8pt,linecolor=blue]{(-}(1.3,-0.8)(1.8,-0.8)
\psline[linewidth=0.8pt,linecolor=blue]{(-|}(1.8,-0.8)(2,-0.8)
\psline[linewidth=0.8pt,linecolor=blue]{(-|}(2.3,-0.8)(2.6,-0.8)
\psline[linewidth=0.8pt,linecolor=blue]{(-}(2.8,-0.8)(3.2,-0.8)
\psline[linewidth=0.8pt,linecolor=blue]{(-|}(3.2,-0.8)(4.3,-0.8)
\psline[linewidth=2pt,linecolor=gray](0,-0.8)(0.1,-0.8)
\psline[linewidth=2pt,linecolor=gray](2,-0.8)(2.3,-0.8)
\psline[linewidth=2pt,linecolor=gray](2.6,-0.8)(2.8,-0.8)
\psline[linewidth=2pt,linecolor=blue](0.1,-0.8)(0.9,-0.8)
\psline[linewidth=2pt,linecolor=blue](0.9,-0.8)(1.3,-0.8)
\psline[linewidth=2pt,linecolor=blue](1.3,-0.8)(1.8,-0.8)
\psline[linewidth=2pt,linecolor=blue]{-}(1.8,-0.8)(2,-0.8)
\psline[linewidth=2pt,linecolor=blue](2.3,-0.8)(2.6,-0.8)
\psline[linewidth=2pt,linecolor=blue](2.8,-0.8)(3.2,-0.8)
\psline[linewidth=2pt,linecolor=blue](3.2,-0.8)(4.3,-0.8)
\rput[b](0.5,-0.75){\small{$1$}}
\rput[b](1.1,-0.75){\small{$2$}}
\rput[b](1.55,-0.75){\small{$3$}}
\rput[b](1.9,-0.75){\small{$4$}}
\rput[b](2.45,-0.75){\small{$5$}}
\rput[b](3,-0.75){\small{$6$}}
\rput[b](3.75,-0.75){\small{$7$}}
\rput[t](0.15,0.2){$\xi_{(1)}^q$}
\rput[t](2.35,0.2){$\xi_{(5)}^q$}
\rput[t](2.85,0.2){$\xi_{(6)}^q$}
\uput[l](4.3,-0.95){\tiny{legend to $F_1^q,\,I_1^q,\,I_2^q\setminus I_1^q,\, I_3^q\setminus I_2^q,\,I_4^q\setminus I_3^q,\,F_2^q,\,I_5^q,\,F_3^q,\,I_6^q$, and $I_7^q\setminus I_6^q$}}
\endpspicture

\vspace{0.1cm}
\centerline{Figure 2}

\smallskip
Suppose that Figure 2 shows the graph of $Z^{\bx,q}$ which corresponds to the same realization as the one illustrated by Figure 1.
The three ``load-free'' intervals $F_i^q$, $i=1,2,3$ are indicated in gray. The interval $I_i^q$ is indicated in blue with marker $i$ on top. 
The excursions of $Z^{\bx,q}$ above past minima are the (closed) disjoint unions of blue intervals. 
Note that $I_i^q$s has exactly the same length as the segment marked by $v(i)$ in Figure 1. Moreover, the following is true.
\begin{Lemma}
\label{Ltosimplify}
Fix a finite initial configuration $\bx$ and time $q>0$. 
If $F_i^q$, $i\geq 1$ are all cut from the abscissa, and the jumps which happen at the end points of $F_\cdot^q$ are all ignored, then (the graph of) $Z^{\bx,q}$ transformed in this way has the law of (the graph of) Aldous' breadth-first walk, corresponding to (\MC) time $q$.
\end{Lemma}
Most of the argument is included in the above made observations and explanations. Figures 1 and 2 illustrate the claim.  The details are left to the reader.

The above reasoning can be summarized as follows:
\begin{Proposition}
\label{Pcoupl}
Let $\bx$ be finite, and fix $q>0$.
Let the breadth-first walk $Z^{\bx,q}$ encode $\bX(q)\in \cvd$ as follows:
for each excursion (above past minima) of $Z^{\bx,q}$ record its length, and let $\bX(q)$ be the vector of thus obtained decreasingly ordered excursion
lengths, appended with infinitely many $0$s. 
Then $\bX(q)$ has the marginal law of the \MC\, started from $\bX(0)=\bx$ and evaluated at time $q$.
\end{Proposition}
Since $(Z^{\bx,q})_{q>0}$ exist on one and the same probability space,  the process 
$$
\bX:=(\bX(q),q>0) \mbox{ and }\bX(0)=\bx
$$
 is well-defined,  and we refer to it temporarily as the {\em \MC\ marginals coupled to $Z^{\bx,\cdot}$}.

{\em Remarks.}
(a)
Apart from simultaneity, the above construction  differs from the one in \cite{aldRGMC,EBMC} in one other significant way:
all the $\xi$s are used ``equally'' in the definition of $Z^{\bx,q}$ for any given $q$ here, while in \cite{aldRGMC,EBMC} for each $q$ separately, the size-biased picking of each leading vertex of a connected component to be explored was done (and encoded) separately. In particular, 
no ``load-free'' intervals exist for the analogues of $Z^{\cdot,\cdot}$ in \cite{aldRGMC,EBMC} (see Figure 1).
At the moment it may seem that the main (potential) gain of the just described modified breadth-first construction is in ``compactifying'' the input data (compare with \cite{EBMC} Section 2.3 or \cite{bromar15}, Section 6.1 for alternatives). It will become apparent in the sequel (see Proposition \ref{Pcoro} and Section \ref{S:SL}) that this construction is quite natural, in that stronger convergence results can be obtained from it with less effort.\\
(b) 
As $q\searrow 0$, the $\xi^q$ diverge to $\infty$, but more importantly they diverge from each other, so $\bX(q) \to  \bx=\bX(0)$ almost surely. 
It is not difficult to see that for any $q\geq 0$, $\bX$ is also almost surely right-continuous at $q$ (see Section \ref{S:BFWUD}).\\
(c)
A little thought is needed to realize that as $q$ increases, the excursion families of $Z^{\bx,\cdot}$ are ``nested'': with probability $1$, if $q_1<q_2$ and two blocks $k,l$ are merged in $\bX(q_1)$, they are also merged in $\bX(q_2)$. This fact is encouraging, but cannot ensure on its own that the \MC\ marginals coupled to $Z^{\bx,\cdot}$ is in fact a \MC. Moreover, while the nesting is encouraging, the following observation will likely increase the level of reader's skepticism about $\bX$ having the \MC\ law: if $e_1^{q_1}$,  $e_2^{q_1}$ and $e_3^{q_1}$ are three different excursions of $Z^{\bx,q_1}$ explored in the increasing order of their indices, and if the initial blocks $k,l,m$ are contained in the connected components matched to $e_1^{q_1},e_2^{q_1},e_3^{q_1}$, respectively, then it is impossible that $k$ and $m$ are merged in $X(q_2)$ without $l$ being merged with $k$ (and therefore with $m$) in $X(q_2)$. 
If there is a simultaneous (for all $q$) scaling limit of $(Z^{\bx,q},\bX(q))$ (under well chosen hypotheses), the just mentioned property persists in the limit.
This observation is perhaps the strongest intuitive argument pointing against the claim of Proposition \ref{Pcoro} and   Theorem \ref{Tmain}.
On the other hand, analogous representations of the standard additive coalescent are well-known (cf.~\cite{chalou,bertoin_frag,bertoin_additive}). 
One may be less surprised there, due to the ``cutting the CRT'' dual (from \cite{dajp97sac}), and the well-known connection between the exploration process of continuum trees and forests on the one hand, and Brownian excursions on the other (cf.~\cite{aldCRT3,pitman-stflour,bertoin-fragcoal}).
As Nicolas Broutin (personal communication) points out,  any (binary) fragmentation can be formally represented as a ``stick-breaking'' process, in which the two broken pieces of any split block remain nearest  neighbors (in some arbitrary but fixed way).
The reversed ``coalescent'' will then have the above counterintuitive property by definition. 
However, one is particularly fortunate if both of these processes (time-reversals of each other) are Markov, and if in addition the ``sticks'' are the excursions of a (generalized) random walk or a related process. \\
(d) Let us denote by ${\mathcal C}$ the operation on paths from Lemma \ref{Ltosimplify}.
In \cite{EBMC} the non-trivial \MC\ extreme entrance laws were obtained by taking limits of Aldous' breadth-first walks (see  Section \ref{S:SL} below), and the limits of their excursions (nearly) above past minima. It was shown that these excursion lengths, considered as an $\cvd$-valued random object, converge in law to the excursion lengths above past minima of the limiting ``walk'' (a member of the family defined in (\ref{defWtc})). 
Lemma \ref{Ltosimplify} makes this latter (somewhat techical) step redundant in the present setting. More precisely, if one can show that under the same hypotheses as those in \cite{EBMC}, $(Z^{t+n^{1/3},\bx^{(n)}})_n$ (see Section \ref{S:SL} for 
precise definitions)
converges to the same $W^{\kappa, t-\tau,\bc}$ as their $({\mathcal C}(Z^{t+n^{1/3},\bx^{(n)}}))_n$, the conclusion about the ordered excursion lengths is immediate. Indeed,  the sequence of excursion above past minima of $Z^{t+n^{1/3},\bx^{(n)}}$
 almost surely matches the sequence of excursions (nearly) above past minima of ${\mathcal C}(Z^{t+n^{1/3},\bx^{(n)}})$, for which Proposition 7 and 9 of \cite{EBMC} (including the results in \cite{EBMC}, Section 2.6) apply verbatim.

\section{Uribe's diagram}
\label{S:UD}
We start by recalling the insight given in Chapter 4 of Uribe \cite{uribe_thesis}, in the notation analogous to that of Section \ref{S:SBFW}. 
In particular, $\bx=(x_1,x_2,\ldots,x_n)$ is a finite-dimensional vector with $n\geq 2$,  $\xi_i$ is an exponential (rate $x_i$) random variable, and $\xi$s are mutually independent.
Denote by $\pi$ the size-biased random reordering of $\bx$, which is determined by $\xi$s, so that $\xi_{\pi_i}\equiv \xi_{(i)}$, $\forall i$ (almost surely). 

Define $n$ different half-lines: for $s\geq 0$
\begin{eqnarray*}
L_1':\bck& &\bck s\mapsto \xi_{(1)} - 0 \cdot s,\\
L_2':\bck& &\bck s \mapsto \xi_{(2)} - x_{\pi_1} s,\\
L_3':\bck& &\bck s \mapsto \xi_{(3)} - (x_{\pi_1} + x_{\pi_2}) s,\\
\bck& &\bck \ldots \ldots \ldots\\
\bck& &\bck \ldots \ldots \ldots\\
L_n':\bck& &\bck s \mapsto \xi_{(n)} - (x_{\pi_1} + x_{\pi_2} + \ldots + x_{\pi_{n-1}}) s,
\end{eqnarray*}
Consider two integers $k,j$ such that $1\leq j<k\leq n$. 
Since $L_k'$ starts (a.s.) at a strictly larger value than $L_j'$, and it has 
 (absolute) slope strictly greater than $L_j'$, it is clear that $L_k'$ and $L_j'$ intersect at some $s_{k,j}>0$.
For each $k=2,\ldots, n$ define
$$
s_k:=\min_{j<k} s_{k,j}, \ \ell_k:=\{1\leq j<k: s_k=s_{k,j}\}.
$$
There are (almost surely) no ties among different $s_k$.
Uribe's diagram consists of line segments (see Figures 3 and 4)
\begin{eqnarray*}
L_1:\bck& &\bck s\mapsto \xi_{(1)} - 0 \cdot s, \ s \in [0, s_2 \vee \ldots \vee s_n],\\
L_2:\bck& &\bck s \mapsto \xi_{(2)} - x_{\pi_1} s, \ s \in [0,s_2],\\
L_3:\bck& &\bck s \mapsto \xi_{(3)} - (x_{\pi_1} + x_{\pi_2}) s, \ s \in [0,s_3], \\
\bck& &\bck \ldots \ldots \ldots\\
\bck& &\bck \ldots \ldots \ldots\\
L_n:\bck& &\bck s \mapsto \xi_{(n)} - (x_{\pi_1} + x_{\pi_2} + \ldots + x_{\pi_{n-1}}) s, \ s \in [0,s_n].
\end{eqnarray*}

\begin{wrapfigure}{r}{0.4\textwidth} 
\vspace{-0.2cm}
    \centering
\includegraphics[width=45mm,height=40mm]{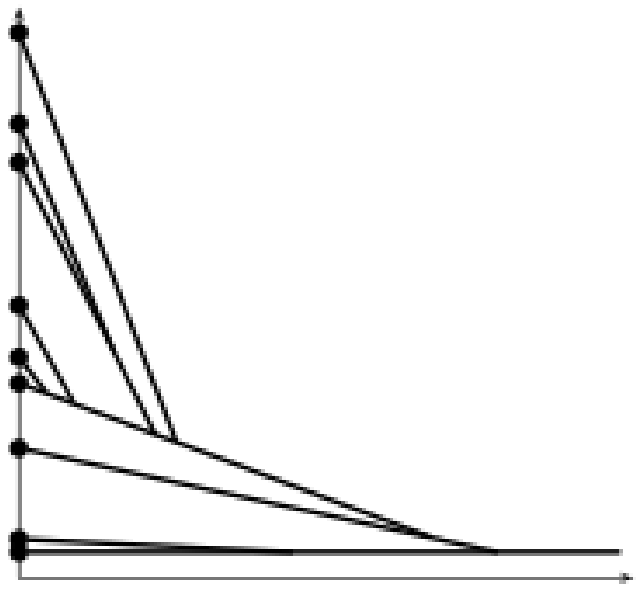}
\centerline{Figure 3}
\vspace{-0.8cm}
\end{wrapfigure}
The image of a realization of $L$ in Figure 3 is inspired by 
 \cite{uribe_thesis}, Chapter 4, Figure 1. 
It is interesting to note that its reflexion in the $x$ axis is \cite{marrat_prep} Figure 1.3.
A more detailed figure, which will correspond in terms of the values of $\bx$, $q$ and $\xi$s to the images in Figures 1 and 2 is provided in Section \ref{S:BFWUD}.

According to \cite{uribe_thesis},  Chapter 4, Section 2,  {\em Armend\'ariz' representation of the \MC} \cite{armen_thesis,armen05} is this graphical summary
 (introduced in \cite{uribe_thesis} for enhanced understanding) joint with an informal description of the following kind:\\
- the state space is $\mathbb{R}_+^n$, its elements are interpreted as lists of block masses, not necessarily ordered,\\
- at time $0$ match the block of mass $x_{\pi_i}$ to $L_i'$,\\
- whenever two of the lines intersect, merge the corresponding blocks, 
form the new vector of block masses accordingly,  continue drawing the lowest indexed line (now matched to the new block), as well as any line that has not participated in the intersection,\\
- keep going as long as there is more than one line remaining.\\
Another construction, advocated as ``essentially Armend\'ariz' representation'' (but better for certain applications) is introduced at the beginning of \cite{uribe_thesis}, Chapter 4, Section 3. This latter construction strongly resembles the ``tilt'' representation of Martin and R\'ath (see \cite{marrat_prep}, Definition 2.7 and Theorem 2.8).

The rest of this section is only one possible rigorous formulation of the aforementioned picture. 
Here Proposition \ref{PUriArm} is stated and proved in a particularly convenient partition-valued framework.
Corollary \ref{Ccoro} features a process that should be the analogue of the Armend\'ariz representation according to \cite{uribe_thesis}. It serves here as a step in obtaining Theorem \ref{Tmain},
through the equivalence obtained in Section \ref{S:BFWUD}, suggesting its potential relevance elsewhere. 
An alternative approach is the particle representation of Martin and R\'ath \cite{marrat_prep}, Section 3.2,  developed in the more general setting of \MC\ with linear deletion. Interestingly, the analysis of \cite{marrat_prep} is again based on an analogue of (a reflection of) Uribe's diagram. 

Uribe's diagram could be interpreted as a ``genealogical tree''. 
More precisely, let us match each point on the diagram $\cup_{i=1}^n\cup_{s\in [0,s_i]} L_i(s)$ to a subset of $\{1,2,\ldots,n\}$ in the following way. For each $i$, set $T_i(0):=\{\pi_i\}$. Each $T_i$ is piece-wise constant, and jumps according to the following algorithm:
$$
T_i(s) :=
 T_i(s-) \,\cup \,
\cup_{j=i+1: s_j=s, \ell_j=i}^n \,T_j(s-), 
 \ i\in \{1,2,\ldots,n\},
$$
and $T_i(s):= \emptyset$, $\forall s \geq s_i$.
In words, for each $i$, the contents of $T_i$ are moved (without replacement) to $T_{\ell_i}$ at time $s_i$, and there is no other copying, cutting or pasting done.
Define $T(s):=\{T_1(s),\ldots, T_n(s)\}$, $s\geq 0$, where the empty sets are ignored.
Note that in this way each $s\geq 0$ is mapped to a partition $T(s)$ of $\{1,2,\ldots,n\}$. 
 The reader is referred to an analogous {\em look-down} construction of \cite{dk99}, which has been since extensively used in the setting of massless (usually called exchangeable) coalescents.

Clearly the partitions along each path of $T$ are nested: if $s_1<s_2$ and $l$ and $k$ are in the same equivalence class of $T(s_1)$, they are (almost surely) in the same  equivalence class of $T(s_2)$.
So $T$ can also be regarded as a random coalescent on the space of partitions of $\{1,2,\ldots, n\}$.
Its initial state  is the trivial partition
$\theta_0:=\{\{\pi_1\},\{\pi_2\},\ldots,\{\pi_n\}\}=\{\{1\},\{2\},\ldots,\{n\}\}$.
Denote by $\GG_t :=\sigma\{T(s), \ s\leq t\}$, $t\geq 0$, so that $\GG:=(\GG_t)_{t\geq 0}$ is the filtration generated by $T$.  The process $T$ will be referred to in the sequel as {\em Uribe's coalescent} process.

It is evident that Uribe's diagram $L$ is a deterministic function of $\bx$ and $\xi$s, and when needed we shall underline this fact by writing $L(\xi_1,\xi_2,\ldots,\xi_n;\bx)$.

{\em Remark.} The $\xi$s are independent and continuous, and therefore (with probability $1$) no two pairs of lines in $L'$ can meet simultaneously. Therefore $T$ (viewed as path-valued) takes value in the space of step functions, such that successive values on a typical path are nested (sub)partitions of $\theta_0$, each having exactly one fewer equivalence class than the prior one. In particular, 
$\GG_t$ is generated by events of the following type: for  $k\geq 1$
$$
\{T(0)=\theta_0,T(t_1)=\theta_1, T(t_2)=\theta_2,\ldots, T(t_k)=\theta_k\}, \ 0< t_1 <\ldots<t_k\leq t,
$$
where $\theta_{j+1}$ is either equal to $\theta_j$ or to a ``coarsening'' of $\theta_j$ obtained by merging two different equivalent classes in $\theta_j$, $0\leq j\leq k-1$. \hfill \hfill\endrem

\smallskip
Let us now account for the masses: for any $i$ and $s\geq 0$ define $M_i(s):=\sum_{l\in T_i(s)} x_l$, with the convention that a sum over an empty set equals $0$. In this way, to each non-trivial equivalence class of $T(s)$ a positive mass is uniquely assigned, and the sum of the masses $\sum_i M_i(s)$ is the identity $\sum_{i=1}^n x_i$, almost surely.

Suppose for a moment that $n=2$. For Uribe's diagram, there are two possibilities: either $\pi$ is the identity, or $\pi$ is the transposition. In either case, the two initial equivalence classes $\{1\}$ and $\{2\}$ merge at random time $s_2$ which we denote by $S$.  Note that the event $\{S>s\}$ is (almost surely) identical to the union of the following two disjoint events $\{\xi_2> \xi_1+s x_1\}$ and $\{\xi_1> \xi_2+s x_2\}$. Thus
\begin{equation}
\label{Etwopart}
P(S>s) = \int_0^\infty x_1e^{-x_1 u} e^{-x_2(u+sx_1)} \,du + \int_0^\infty x_2e^{-x_2 u} e^{-x_1(u+sx_2)} \,du,
\end{equation}
and the reader can easily verify that the RHS equals $e^{-x_1x_2s}$. So if $n=2$, the coalescent time is distributed equally in the random graph and in Uribe's coalescent. 

Even with this hint in mind, the next result will likely seem at least counterintuitive if not  striking to an uninitiated reader.
\begin{Proposition}
\label{PUriArm}
Uribe's coalescent process $T$ has the random graph law. More precisely, it is a continuous-time Markov chain, such that any two equivalence classes in $T$ merge, independently of all the other merger events, at the rate equal to the product of their masses.
\end{Proposition}
The argument given below is based on a sequence of elementary observations, and
its outline is comparable to that of \cite{armen_thesis} Lemma 15.
Some care is however needed in correctly setting up the conditioning (otherwise the statement of the proposition would seem obvious from the start). In particular, there seems to be no way of a priori knowing that $T$ with respect to $\GG$ has the Markov property. The argument included below exhibits the transition rates in the process of checking for Markovianity. 
 
\medskip
{\em Proof of Proposition \ref{PUriArm}.}
The random graph on finitely many vertices evolves by exponential jumps, and it changes states whenever, due to an arrival of a new edge between two original blocks, the connectivity of the graph changes (in that one new connected component is formed from two previous ones). The new edges that appear but do not change the connectivity can (and will be) ignored for our purposes.

Now consider $n\geq 3$. The minimal coalescence time in the random graph has exponential distribution with rate
$\sum_{i,j: 1\leq i< j\leq n} x_i x_j$.
Let $S:=\min\{s_2,s_3,\ldots,s_n\}$ be the corresponding minimal coalescence time in Uribe's coalescent. One can extend (\ref{Etwopart}) by noting that 
$\{S>s\}=\{T(s)=\theta_0\}$ $=\{\xi_{\pi_2}>\xi_{\pi_1}+sx_{\pi_1},\,\xi_{\pi_3}>\xi_{\pi_2}+s x_{\pi_2},\,\ldots,\xi_{\pi_n} >\xi_{\pi_{n-1}} + s x_{\pi_{n-1}}\}$ can be split into $n!$ disjoint events according to the value of the random permutation $\pi$. Fix $\tau$ a deterministic permutation and note that 
$P(S>s, \pi=\tau)$ equals, similarly to (\ref{Etwopart}),
\begin{eqnarray}
\label{Enpart}
& &\bck\bck\int_0^\infty \bck du_1\,x_{\tau_1}e^{-x_{\tau_1} u_1} \int_{u_1+sx_{\tau_1}}^\infty \bck du_2\, x_{\tau_2} e^{-x_{\tau_2} u_2}
\int_{u_2+sx_{\tau_2}}^\infty \bck du_3\, x_{\tau_3} e^{-x_{\tau_3} u_3} \ \cdots\nonumber\\
 & &\ \cdots\int_{u_{n-2}+sx_{\tau_{n-2}}}^\infty \bck du_{n-1}\, x_{\tau_{n-1}} e^{-x_{\tau_{n-1}} u_{n-1}}  \int_{u_{n-1}+sx_{\tau_{n-1}}}^\infty \bck du_n\, x_{\tau_n} e^{-x_{\tau_n} u_n} =\nonumber\\
& &\bck\bck e^{-x_{\tau_{n-1}} x_{\tau_n} s}\int_0^\infty \bck du_1\,x_{\tau_1}e^{-x_{\tau_1} u_1} \int_{u_1+sx_{\tau_1}}^\infty \bck du_2\, x_{\tau_2}  e^{-x_{\tau_2} u_2} \int_{u_2+sx_{\tau_2}}^\infty\bck  du_3 x_{\tau_3}  e^{-x_{\tau_3} u_3}\nonumber\\
& & 
 \ \cdots \int_{u_{n-2}+sx_{\tau_{n-2}}}^\infty du_{n-1}\, x_{\tau_{n-1}} e^{-(x_{\tau_{n-1}} + x_{\tau_n}) u_{n-1}} =\nonumber\\
& &\bck\bck e^{-x_{\tau_{n-1}} x_{\tau_n} s} \frac{x_{\tau_{n-1}}}{x_{\tau_{n-1}} + x_{\tau_n}} e^{-x_{\tau_{n-2}}(x_{\tau_{n-1}} + x_{\tau_n}) s}
\int_0^\infty \bck du_1\,x_{\tau_1}e^{-x_{\tau_1} u_1} \cdot  \nonumber
\end{eqnarray}
\begin{eqnarray}
& & \int_{u_1+sx_{\tau_1}}^\infty \bck du_2\, x_{\tau_2}  e^{-x_{\tau_2} u_2} \int_{u_2+sx_{\tau_2}}^\infty  \bck du_3\, x_{\tau_3} e^{-x_{\tau_3} u_3}\cdots \nonumber \\
& & \ \cdots \int_{u_{n-3}+sx_{\tau_{n-3}}}^\infty du_{n-2}\, x_{\tau_{n-2}} e^{-(x_{\tau_{n-2}} + x_{\tau_{n-1}} + x_{\tau_n}) u_{n-2}} = \cdots
\label{Emanyend}
\end{eqnarray}
We proceed  inductively with integration in (\ref{Emanyend}) to obtain that  $P(S>s, \pi=\tau)$ equals
$$
\prod_{i=1}^{n-1} \frac{x_{\tau_i}}{\sum_{j=i}^n x_{\tau_j} } e^{-s x_{\tau_i}\cdot \sum_{j=i+1}^n x_{\tau_j}} = 
\prod_{i=1}^{n-1} e^{-s x_i\cdot \sum_{j=i+1}^n x_j} \prod_{i=1}^{n-1} \frac{x_{\tau_i}}{\sum_{j=i}^n x_{\tau_j}},
$$
and this final quantity can be recognized as 
\begin{equation}
\label{ESpilaw}
e^{- s\sum_{i,j: 1\leq i< j\leq n} x_i x_j} P(\pi= \tau).
\end{equation}
An immediate conclusion is that 
\begin{equation}
\label{EindepSpi}
\begin{array}{c}
S \mbox{ is independent of }\pi \mbox{ in Uribe's coalescent,  and}\\
\mbox{$S$ has the same law as its counterpart in the random graph.}
\end{array}
\end{equation}
The just obtained equivalence is a good start, but we need a stronger form of independence in order to obtain the full equivalence in law.
In view of this, let us consider the family of residuals $(\xi_i - s\sum_{j:\xi_j<\xi_i} x_j,\,i=1,2,\ldots,n)$ on the event 
$\{S>s\}$. Note that $L_k(s)$ equals the $i$th residual time above if and only if $\pi_k=i$ or $k=\pi^{-1}(i)$. It is natural to denote $L_{\pi^{-1}(i)}(s)$ by $\xi_i(s)$. The calculation leading to (\ref{ESpilaw}) can be only slightly modified to yield that for $v_i>0$, $i=1,2,\ldots,n$ we have
$P(S>s, \pi=\tau, \xi_i(s) = dv_i,\,i=1,2,\ldots,n) =$ $ P(\xi_i = dv_i+ s \sum_{j:\tau^{-1}(j)<\tau^{-1}(i)}x_j,\,i=1,2,\ldots,n) =$
$$
\prod_{i=1}^n \left(x_i e^{-x_i v_i}\,dv_i \cdot e^{- s\,x_i\cdot\sum_{j: \tau^{-1}(j)<\tau^{-1}(i)} x_j}\right) \cdot P( \pi=\tau).
$$
But the product of the exponential terms is again $P(S>s)$, so we see that 
given $\GG_s$ and on $\{S>s\}$,
the residual random variables $\{\xi_i(s),\,i=1,\ldots,n\}$ have the conditional joint distribution equal to the original distribution of $\{\xi_i(0)=\xi_i,\,i=1,\ldots,n\}$. The order of indices induced by $\xi(s)$s is the same as that induced by $\xi$s.
So given $\GG_s$, and knowing $\{S>s\}=\{T(s)=\theta_0\}$, the residual Uribe's diagram $(L(s+u),\, u\geq 0)$, where $L=L(\xi_1(0),\xi_2(0),\ldots,\xi_n(0);\bx)$, has the  conditional law equal to that of Uribe's diagram $L(\xi_1(s),\xi_2(s),$ $\ldots,\xi_n(s);\bx)$, which in this case equals the law of $L$ (as argued already, on $\{T(s)=\theta_0\}$ the $(\xi_i(s))_i$ are again independent exponentials with respective rates $(x_i)_i$). Therefore,
\begin{equation}
\label{Econdlawequiv}
\begin{array}{c}
 \mbox{ given }\GG_s \mbox{ and on }\{T(s)=\theta_0\}, \mbox{ the process }(T(s+u),\, u\geq 0)\\
\mbox{ has the conditional law equal to }(T(u),\,u\geq 0).
\end{array}
\end{equation}

The proof will be completed if we can verify an analogous statement for any possible event in $\GG_s$. 
In fact it suffices to concentrate on the class of  events from the remark stated after the definition of $\GG$, and show that the laws of Uribe's coalescent started from $\theta_0$, and the random graph started from $n$ blocks with respective sizes $x_1,\ldots, x_n$, agree on each event of the form given in that remark.

The next few paragraphs serve to argue an intermediate step, accounting for a single merger recorded up to the present time. One should show that, in this special case, the future depends on the past only through the present.
Without loss of generality we can pretend (or reindex the blocks so) that the first jump of $T$ is from $\theta_0$ to $\theta_{12}:=\{\{1,2\},\{3\},\ldots,\{n\}\}.$
Let ${\mathcal S}_{12}$ be the subset of permutations of $\{1,2,\ldots,n\}$ defined by
$$\tau\in {\mathcal S}_{12}\mbox{ iff }|\tau^{-1}(1)-\tau^{-1}(2)|=1.$$
Consider the event $\{T(s)=\theta_{12}\}=\{T(0)=\theta_0,T(s)=\theta_{12}\}$. 
In the notation of the remark following the definition of $\GG$, one has $t_1=s=t$ and $\theta_1=\theta_{12}$.
It is clear that  $\{T(s)=\theta_{12}\}\subset \{\pi \in {\mathcal S}_{12}\}$, or equivalently, that
$\{T(s)=\theta_{12}\}$ $=\cup_{\tau\in  {\mathcal S}_{12}} 
\{T(s)=\theta_{12}, \pi=\tau\}$. 

For $\tau\in {\mathcal S}_{12}$, 
denote by $\tau^*$ a natural transposition of $\tau$: $\tau^*_j=\tau_j$, $\forall j$ such that $\tau_j\not\in\{1,2\}$, and  $\tau^*_j=3-\tau_j$, if $\tau_j\in \{1,2\}$.
Now suppose that $\tau_i=1=3-\tau_{i+1}$ for some $i<n-1$ (the upper limit $n-1$ for $i$ is treated separately).
Then on $\{T(s)=\theta_{12}, \pi=\tau\}$ (resp.~$\{T(s)=\theta_{12}, \pi=\tau^*\}$) there must exist  $u\leq s$ such that $\xi_2-\xi_1= u x_1$ (resp.~$\xi_1-\xi_2= u x_2$) and in addition it must be true that $\xi_{\tau_{j+1}} > \xi_{\tau_j} + s x_{\tau_j}$, for all $j=1,\ldots, i-1, i+2,\ldots,n$ as well as $\xi_{\tau_{i+2}} - \xi_1 \,\wedge\,  \xi_2> (x_{\tau_i} + x_{\tau_{i+1}})s= (x_{\tau^*_i} + x_{\tau^*_{i+1}})s= (x_1+x_2) s$.
So, in analogy to (\ref{Enpart}), $P(T(s)=\theta_{12},\pi=\tau)$ becomes
\begin{eqnarray}
& &\bck\bck\int_0^\infty \bck du_1\,x_{\tau_1}e^{-x_{\tau_1} u_1} \int_{u_1+sx_{\tau_1}}^\infty \bck du_2\, x_{\tau_2} e^{-x_{\tau_2} u_2}
 \cdots\nonumber\\
& &\bck \bck\cdot\int_{u_{i-1}+sx_{\tau_{i-1}}}^\infty \bck du_i\, x_1\, e^{-x_1 u_i}  \int_{u_i}^{u_i+sx_1} \bck du_{i+1}\, x_2\, e^{-x_2 u_{i+1}} \cdot\nonumber\\
& &\cdot \int_{u_i+s(x_1+x_2)}^\infty \bck x_{\tau_{i+2}} e^{-x_{\tau_{i+2}} u_{i+2}} \cdots \int_{u_{n-1}+sx_{\tau_{n-1}}}^\infty \bck du_n \, x_{\tau_n} e^{-x_{\tau_n} u_n}, \label{Enpart a}
\end{eqnarray}
and the symmetric quantity for $P(T(s)=\theta_{12},\pi=\tau^*)$ differs from the above expression only in the middle row above as follows:
$$
\bck \bck\cdot\int_{u_{i-1}+sx_{\tau_{i-1}}}^\infty \bck du_i\, x_2\, e^{-x_2 u_i}  \int_{u_i}^{u_i+sx_2} \bck du_{i+1}\, x_1\, e^{-x_1 u_{i+1}} \cdot
$$
Adding the two integral expressions, and evaluating the integral, we get that $P(T(s)=\theta_{12},\pi\in\{\tau,\tau^*\})$ equals
\begin{equation}
\label{Econcl}
(1-e^{-sx_1x_2}) \cdot e^{- s\sum_{j,k: 1\leq j< k\leq n, (j,k)\neq (1,2) } x_j x_k} \cdot I_{12}^{\tau,\tau^*}
\end{equation}
where 
\begin{equation}
\label{EI12}
 I_{12}^{\tau,\tau^*}=  I_{12}^{\tau^*,\tau}=
\prod_{j=1}^{i-1} \frac{x_{\tau_j}}{\sum_{k=j}^n x_{\tau_j}}\cdot 
\frac{x_1+x_2}{x_1+x_2 + \sum_{j=i+2}^n x_{\tau_j}} \cdot 
\prod_{j=i+2}^{n-1} \frac{x_{\tau_j}}{\sum_{k=j}^n x_{\tau_j}}.
\end{equation}
The assumption that $i+1<n$ implies the presence of the final integral in the second line, and possibly extra terms in the third line of (\ref{Enpart a}). In the special case $i=n-1$, the analogous formulae are somewhat simpler, but the conclusion (\ref{Econcl}) is the same (as the reader can easily verify) with $I_{12}^{\tau.\tau^*}$ redefined as 
$\prod_{j=1}^{n-2} \frac{x_{\tau_j}}{\sum_{k=j}^n x_{\tau_j}}$. Since the sum of $I_{1.2}$s over all the different choices of permulation (and its transposition)  $\tau,\tau^*\in {\mathcal S}_{12}$ is clearly equal to $1$, one concludes moreover that 
\begin{equation}
\label{Edistpiprime}
P(\pi\in\{\tau,\tau^*\}|T(s) =\theta_{12})=  I_{12}^{\tau,\tau^*},
\end{equation}
and that $P(T(s)=\theta_{12})=(1-e^{-sx_1x_2}) \,\exp\{- s\sum_{j,k: 1\leq j< k\leq n, (j,k)\neq (1,2) } x_j x_k\}$. 

Define a map 
$$' \, : \, {\mathcal S}_{12} \mapsto\mbox{ permutations of }\{12,3,\ldots,n\}$$
as follows: for $k\in \{1,2,\ldots,n-1\}$ let
\[
\tau'(k):=\left\{ 
\begin{array}{ll}
\tau(k), & k< \min\{\tau^{-1}(1),\tau^{-1}(2)\}, \\
12, & k =\min\{\tau^{-1}(1),\tau^{-1}(2)\},\\
\tau(k+1), & k \geq \min\{\tau^{-1}(1),\tau^{-1}(2)\}+1.
\end{array}
\right.
\]
In words, $\tau'$ is obtained from $\tau\in {\mathcal S}_{12}$ by checking which of the two indices 1 or 2 occurs first in $\tau$, in renaming this image point $12$, in deleting the next image point (which is necessarily 1 or 2), and in  preserving the image order set by $\tau$ otherwise. 
The identity in (\ref{Econcl}) can now be restated as 
\begin{equation}
\label{EindepSLpiprime}
\begin{array}{c} 
\{T(s)=\theta_{12}\}   \mbox{ has equal probability in the random graph and in Uribe's coalescent,}
\\
\mbox{ and, on $\{T(s)= \theta_{12}\}$, the conditional law of $\pi'$ given $\GG_s$ is specified by (\ref{EI12},\ref{Edistpiprime}).}
\end{array}
\end{equation}
In fact more is true, in analogy to the behavior already verified with respect to $\{T(s)=\theta_0\}$. 
Define $\xi_{12}:=\xi_1 \,\wedge\,\xi_2$ and consider the family of $n-1$ residuals $\{(\xi_i(s):=\xi_i - s\sum_{j:\xi_j<\xi_i} x_j,\,i=3,\ldots,n),\xi_{12}(s):=\xi_{12}- s \sum_{j:\xi_j<\xi_{12}} x_j\}$.
On $\{T(s)=\theta_{12}\}$ the permutation of indices $\{12,3,\ldots,n\}$ induced by $\xi_{12}(s)$, $\xi_3(s),\ldots$, $\xi_n(s)$ is precisely $\pi'$.
Note that $\xi_j< \xi_{12}$ if and only if $\pi^{-1}(j)<\pi^{-1}(1)\wedge \pi^{-1}(2)$, or if and only if $\pi'^{-1}(j)<\pi'^{-1}(12)$.
Now let $\tau \in {\mathcal S}_{12}$ be such that $\tau^{-1}(1)=\tau^{-1}(2)-1$.
Then $P(T(s)=\theta_{12}, \pi=\tau, \xi_i(s) = dv_i,\,i\in\{12,3,\ldots,n\}) =$ 
$ P(\xi_1=dv_{12} + s \sum_{j:\tau^{-1}(j)<\tau^{-1}(1)}x_j, \xi_i = dv_i+ s \sum_{j:\tau^{-1}(j)<\tau^{-1}(i)}x_j,\,i\in\{3,\ldots,n\}, \,\xi_2- \xi_1 \leq s\, x_1)$, and this
equals
\begin{eqnarray}
x_1 e^{-x_1 v_{12}} e^{ - x_1\cdot s \sum_{j:\tau^{-1}(j)<\tau^{-1}(1)}x_j} \cdot  
\prod_{i=3}^n \left(x_i e^{-x_i v_i} \cdot 
e^{- x_i\cdot s\sum_{j: \tau^{-1}(j)<\tau^{-1}(i)}x_j}\right)\cdot\nonumber\\
\int_{v_{12}+  s \sum_{j:\tau^{-1}(j)<\tau^{-1}(1)}x_j}^{v_{12}+  s \left(\sum_{j:\tau^{-1}(j)<\tau^{-1}(1)}x_j +x_1\right)}x_2e^{-x_2 u} du\, dv_{12}\, dv_3\cdots dv_n=\nonumber \\
x_1 e^{-(x_1+x_2) v_{12}} \prod_{i=3}^n x_i e^{-x_i v_i}\,dv_{12}\, dv_3\cdots dv_n \, \cdot (1-e^{-sx_1x_2}) \cdot  \nonumber \\  
\cdot \exp\left\{- s\left(\left[\sum_{i=3}^n \sum_{j: \tau^{-1}(j)<\tau^{-1}(i)}x_i x_j\right] + (x_1+x_2) \cdot \bck\bck\bck\sum_{j: \tau^{-1}(j)<\tau^{-1}(1)} x_j\right)\right\}.\ 
\label{Eterms12}
\end{eqnarray}
The final exponential in (\ref{Eterms12}) is easily seen to be equal to the exponential in (\ref{Econcl}).
An analogous expression can be written for $P(T(s)=\theta_{12}, \pi=\tau^*, \xi_i(s) = dv_i,\,i\in\{12,3,\ldots,n\})$, and by symmetry considerations, the reader can easily see that the only difference between the two products is the leading multiple which changes from $x_1$ to $x_2$.  The conclusion is that 
$P(T(s)=\theta_{12}, \pi'=\tau', \xi_i(s) = dv_i,\,i\in\{12,3,\ldots,n\})$ equals
$$
(x_1+x_2) e^{-(x_1+x_2) v_{12}} \prod_{i=3}^n x_i e^{-x_i v_i}\,dv_{12}\, dv_3\cdots dv_n \cdot P(T(s)=\theta_{12}),
$$
and that therefore, given $\GG_s$, and knowing $\{T(s)=\theta_{12}\}$, the residual Uribe's diagram
$(L(s+u),\, u\geq 0)$, where $L=L(\xi_1(0),\xi_2(0),\ldots,\xi_n(0);\bx)$, has the conditional law of $L(\xi_{12}(s),\xi_3(s),\ldots,\xi_n(s);(x_1+x_2,x_3,\ldots,x_n))$, which is Uribe's diagram generated from $n-1$-dimensional vector 
$(x_1+x_2,x_3,\ldots,x_n)$ and its corresponding independent exponentials $\xi_\cdot(s)$s.

Using (\ref{EindepSpi}), (\ref{Econdlawequiv}), (\ref{EindepSLpiprime}), and the just made conclusion, 
a standard (inductive) nested conditioning argument yields that for any $k\in \mathbb{N}$, any $0< t_1 <\ldots<t_k<\infty$, and any sequence $(\theta_j)_j$ where,  for each $j\leq k-1$, $\theta_{j+1}$ either equals $\theta_j$ or a single merger coarsening of $\theta_j$, it is necessarily true that the event
$$
\{T(0)=\theta_0,T(t_1)=\theta_1, T(t_2)=\theta_2,\ldots, T(t_k)=\theta_k\}
$$
occurs with equal probability in Uribe's coalescent and in the random graph. 
As already argued, this is sufficient to conclude the full identity in law. \hfill \endpf

Recall that, for $i=1,2,\ldots,n$, $M_i(s)$ is the mass of the equivalence class $T_i(s)$, provided that $T_i(s)\neq \emptyset$, and it is defined to be $0$ otherwise. 
Now let $\bY(s)$ be a $\cvd$-valued random variable, formed by listing the components of $(M_1(s),M_2(s),\ldots,M_n(s))$ in decreasing order, and appending infinitely many zeros (to obtain a vector in $\cvd$). 
It is clear that $\bY=(\bY(s),\,s\geq 0)$ is adapted to $\GG$. Moreover Proposition \ref{PUriArm} can be restated as 
\begin{Corollary}
\label{Ccoro}
The process $\bY$ is a \MC\ started from the decreasing ordering of $(x_1,x_2,\ldots,x_n,0,\ldots)$.
\end{Corollary}

In the case where all the $n$ initial masses are equal, a subset of the just derived identities was already known to Gumbel \cite{gumbel}. 
It is well-known (see for example the discussion in \cite{uribe_thesis}, Chapter 4, or \cite{bol_book}, Chapter 7, or \cite{durrett_rgd}, Chapter 2, Section 8) that the connectivity time of the (classical) Erd\"os-R\'enyi random graph, is of the order $(\log{n}+G+o(1))/n$, where $G$ has Gumbel's law $P(G\leq g) = e^{-e^{-g}}$, $g\in\mathbb{R}$.

\section{Breadth-first walks meet Uribe's diagram}
\label{S:BFWUD}
In this section we will compare the simultaneous breadth-first random walks of Section \ref{S:SBFW} with Uribe's diagram of Section \ref{S:UD}. More precisely, a coupling of these random objects will be realized on one and the same probability space, so that the \MC\ marginals $\bX$ coupled to $Z^{\bx,\cdot}$ (see Section \ref{S:SBFW}) can be matched to $\bY$ derived from Uribe's diagram (see Section \ref{S:UD}).  

As an immediate corollary we obtain 
\begin{Proposition}
\label{Pcoro}
Let $\bx$ be finite. Then the \MC\ marginals $\bX$ coupled to $Z^{\bx,\cdot}$ has the law of a \MC\ started from $\bx$. 
\end{Proposition}

Let $\bx$ be finite, set $n:=\len(\bx)$, and recall the construction from Section \ref{S:SBFW}. Use the same $(\xi_1,\ldots,\xi_n)$ to form the corresponding Uribe diagram $L(\xi_1,\xi_2,\ldots,\xi_n;\bx)$. The notation is slightly abused here since in Section \ref{S:SBFW} (resp.~\ref{S:UD}) vectors have infinite (resp.~finite) length, but the correspondence between the two is clear, whenever there are only finitely many blocks of non-trivial mass in the initial state. 

Assume for a moment that $\bx=(1.1,0.8,0.5,0.4,0.4,0.3,0.2)$, and let us pretend that 
the realization from Figures 1 and 2 corresponds to 
the time parameter $q$ equal to $2$. 
Let us assume in addition that $\pi=\tau:=(2,4,3,7,6,5,1)$.
This means that $(\xi_{(1)},\xi_{(2)},\ldots, \xi_{(7)})=(\xi_2,\xi_4,\xi_3,\xi_7,\xi_6,\xi_5,\xi_1)$, where the latter vector takes value $(0.2, 0.7, 1.4, 3.4, 4.6, 5.6, 6)$.  

The corresponding Uribe's diagram $L(\xi,\bx)$ is shown in Figure 4 below.
For any time $s$, the partition $T(s)$ can be read from the graph as it could be read from a genealogical tree.
Each of the ``active'' lines represents a different equivalence class.
The blue vertical dashed line marks time $s=2$. 
In $T(2)$ there are three equivalence classes, matched to $L_1, L_5$ and $L_6$, as shown in the figure. 
This partition is, of course, the same as the one depicted on Figure 1 (recall that the same realization is being illustrated in Figures 1, 2 and 4).

The coupling stated at the beginning of the section is realized in the most natural way.
Recall that both $(Z^{\bx,q}, \, q>0)$ and $L(\xi;\bx)$ (and therefore $T$ and $\bY$) are functions of $\xi$s and $\bx$.
Let the finite family of independent exponential random variables $\xi_\cdot$, used in the construction of $Z^{\bx,\cdot}$ and $L(\xi;\bx)$, be the same, almost surely.

  \psset{linewidth=0.4pt}
  \psset{xunit=2.762cm,yunit=2.9cm}
   \pspicture(-0.2,-0.5)(4.4,6.4)
 \psset{linestyle=solid}
   \psline[linecolor=gray]{-}(0.625,0.2)(1,-0.1)
 \psline[linecolor=gray]{-}(1,0.2)(1.25,-0.1)
\psline[linecolor=gray](1.88235,0.2)(2.1,-0.17)
\psline[linecolor=gray](2.31578,0.2)(2.5,-0.15)
  \psline[linecolor=gray]{*-}(2.4545,0.2)(2.65,-0.23)  
 \psline[linecolor=gray]{-}(1,3.4)(1.5,2.1)              
   \psset{linestyle=solid}
   \psline{->}(0,6.3)\psline{-}(0,-0.3)
  \psline{->}(4.5,0)
   \uput[l](4.35,-0.07){\small{$s$}}
    \psset{linewidth=1pt}
    \psline{*-}(0,0.2)(4,0.2)            
    \psline{*-*}(0,0.7)(0.625,0.2)          
    \psline{*-*}(0,1.4)(1,0.2)          
    \psline{*-*}(0,3.4)(1.88235,0.2)          
    \psline{*-*}(0,4.6)(2.31578,0.2)          
    \psline{*-*}(0,5.6)(2.4545,0.2)          
    \psline{*-*}(0,6)(1,3.4)             
\psline[linestyle=dashed,linecolor=blue](2,-0.2)(2,2)
\pnode(2,1.2){third}
\pnode(2,0.8){second}
\pnode(2,0.2){first}
\rput[b](3,0.3){\rnode{A}{\small{$T_1(2)=\{2,4,3,7\}$}}}
\rput[b](2.7,0.8){\rnode{B}{\small{$T_5(2)=\{6\}$}}}
\rput[b](2.7,1.2){\rnode{C}{\small{$T_6(2)=\{5,1\}$}}}
\ncarc[linewidth=0.4pt,nodesepA=2pt,nodesepB=3pt,arcangleA=-90, arcangleB=-45]{->}{A}{first}
\ncarc[linewidth=0.4pt,nodesepA=4pt,nodesepB=3pt,arcangleA=145, arcangleB=35]{->}{B}{second}
\ncarc[linewidth=0.4pt,nodesepA=4pt,nodesepB=3pt,arcangleA=135,arcangleB=35]{->}{C}{third}
\rput[r](-0.1,0.2){\small{$\{2\}$}}
\rput[r](-0.1,0.7){\small{$\{4\}$}}
\rput[r](-0.1,1.4){\small{$\{3\}$}}
\rput[r](-0.1,3.4){\small{$\{7\}$}}
\rput[r](-0.1,4.6){\small{$\{6\}$}}
\rput[r](-0.1,5.6){\small{$\{5\}$}}
\rput[r](-0.1,6){\small{$\{1\}$}}
\endpspicture

\vspace{-0.5cm}
\centerline{Figure 4}

\medskip
As already noted (see Remark (c) at the end of Section \ref{S:SBFW}), the partition structure induced by the evolution of $Z^{\bx,\cdot}$ gets coarser as $q$ increases. In addition, only pairs of neighboring blocks or families of blocks, with respect to the random order established by $\pi$, can coalesce either in $\bX$ or in $T$ (that is, in $\bY$). 
Note that, for each \MC\ time $q$, the relation of being connected by a path of edges $ \lrarrow$ that occurred before time $q$ is an equivalence relation on the initial set of blocks.
Hence it suffices to show that, almost surely, for each $q>0$ and $i\in \{1,2,\ldots,n-1\}$ it is true that 
$\pi_{i} \lrarrow \pi_j$ with respect to $Z^{\bx,q}$ (see Section \ref{S:SBFW}) if and only if $\pi_i \sim \pi_j$ with respect to $T$ (see Section \ref{S:UD}).
At time $q=0$ the just made claim is clearly correct, since there are no edges $\lrarrow$ in $\bX$, and the partition of $T$ is trivial.
Suppose that random time $Q_1>0$ is such that $T(Q_1-)= \theta_0$ and $T(Q_1)$ contains  $\{\pi_i,\pi_{i+1}\}$. This means that  the lines $L_i'$ and $L_{i+1}'$ intersect at time $Q_1$, and no other pair of lines intersects before time $Q_1$.
Or equivalently, $\xi_{(i+1)}-  \xi_{(i)}= Q_1 x_{\pi_i}$ and $\xi_{(j+1)}-  \xi_{(j)}> Q_1 x_{\pi_j}$ for $j\neq i$.
Or equivalently, $\xi_{(i+1)}^{Q_1} -  \xi_{(i)}^{Q_1}= x_{\pi_i}$, and $\xi_{(j+1)}^{Q_1} -  \xi_{(j)}^{Q_1}>x_{\pi_j}$ for  $j\neq i$.
A quick check of the construction in Section \ref{S:SBFW} suffices to see that, on the above event, the edge  $\pi_i \lrarrow \pi_{i+1}$ arrives at time $Q_1$ in $\bX$, and no edge arrives to $\bX$ before time $Q_1$.
At  time $Q_1$, the line $L_{i+1}'$ stops being active in Uribe's diagram. The new neighbors of $\{\pi_i,\pi_{i+1}\}$ are $\pi_{i+2}$ and $\pi_{i-1}$. All the active lines above  $L_i'$ account for the new mass $x_{\pi_i}+x_{\pi_{i+1}}$ of  $\{\pi_i,\pi_{i+1}\}$, since this quantity is built into their slope (together with masses corresponding to any other active lines underneath them). 
Similarly,  $Z^{\bx,q}$ for $q>Q_1$ does not need to observe $\xi_{(i+1)}^q$ any longer, it suffices to attribute the cumulative ``listening length'' $x_{\pi_i}+x_{\pi_{i+1}}$ to the breadth-first walk time $\xi_{(i)}^q$ at which the leading particle of the component $\{\pi_i,\pi_{i+1}\}$ is seen by the walk.
Due to these two observations, one can continue the comparison of the coalescence of the remaining blocks driven by $(T(q),\,q\geq Q_1)$ to that driven by $(Z^{\bx,q},\,q\geq Q_1)$, and conclude by induction that in both processes the sequence of pairs of blocks that coalesce, and their respective times of coalescence, are identical, almost surely.\\
As already noted, Proposition \ref{Pcoro} is a direct consequence. 
Note that in the sense of the just produced coupling, 
the simultaneous breadth-first walks of Section \ref{S:SBFW} are equivalent to Uribe's diagram. 
In comparison, the original  breadth-first walk of \cite{aldRGMC} coupling is ``static'' (it works for one $q$ at a time), and it seems difficult to turn it into a ``dynamic'' version due to a certain (small but present) loss of information (see Lemma \ref{Ltosimplify}). 

\section{Scaling limits for simultaneous breadth-first walks}
\label{S:SL}
This section imitates the approach of Section 2.4 in \cite{EBMC}. 
It is interesting to note that these scaling limits are simpler to derive here than they were for the original \MC\ encoding walks in \cite{EBMC}.

Given $\bx \in \cvd$ let
\[ \sigma_r(\bx) := \sum_i x_i^r, \ r = 1,2,3. \]
For each $n \geq 1$ let $\bx^{(n)}$ be a finite vector (in the sense of Section \ref{S:SBFW}).
Let $((Z^{\bx^{(n)},q}(s), \, s\geq 0),\,q\geq 0)$ and $(\bX^{(n)}(q), q \geq 0)$ be the simultaneous breadth-first walks, and  
 the \MC\ coupled to  $Z^{\bx^{(n)},\cdot}$ (see Section \ref{S:SBFW} and Proposition \ref{Pcoro}), respectively. 

Suppose that for some $\kappa \in[0,\infty)$ and  $\bc \in \cvt$, the following hypotheses are true:
\begin{eqnarray}
\frac{\sigma_3(\bx^{(n)})}{(\sigma_2(\bx^{(n)}))^3} &\con& \kappa + \sum_j c_j^3, \label{hyp1}\\
 \frac{x_j^{(n)}}{\sigma_2(\bx^{(n)})} &\con& c_j, \ j \geq 1,  \label{hyp3}\\
\sigma_2(\bx^{(n)}) &\con& 0 ,\label{hyp4}
\end{eqnarray}
 as $n \con \infty$.
It is easy to convince oneself (or see  Lemma 8 of \cite{EBMC}) that for any
$(\kappa,0,\bc) \in {\cal I}$
there exists a finite vector valued sequence $(\bx^{(n)})_{n\geq 1}$ satisfying (\ref{hyp1} -- \ref{hyp4}).
In the standard Aldous' setting of \cite{aldRGMC}, $x^{(n)}$ had length $n$, and all of its blocks had identical mass $1/n^{2/3}$. 
If $\bc\neq {\bf 0}$,  it is then natural to introduce ``large'' dust blocks of mass of order $1/n^{1/3}$.

As in \cite{EBMC}, we furthermore pick  an integer valued sequence $(m(n))_{n\geq 1}$, which increases to infinity 
sufficiently slowly so that
\begin{equation}
\label{mChoo1new}
\left| \sum_{i=1}^{m(n)} \frac{(x_i^{(n)})^2}{(\sigma_2(\bx^{(n)}))^2}
  - \sum_{i=1}^{m(n)} c_i^2 \right| \con 0 \;,\; \; 
\left| \sum_{i=1}^{m(n)} \left( \frac{x_i^{(n)}}{\sigma_2(\bx^{(n)})}
  - c_i \right)^3 \right| \con 0 \;,\; \; 
\end{equation}
\begin{equation}
\label{mChoo2new}
\mbox{ and } \; \sigma_2(\bx^{(n)})\sum_{i=1}^{m(n)} c_i^2 \con 0 \;.
\mbox{ \hspace{1cm} }
\end{equation}

Fix $t\in \mathbb{R}$ and let $q_n:= \sfrac{1}{\sigma_2(\bx^{(n)})} + t$. 
Recall (\ref{defZbxq}), and define
$$Z_n:=Z^{\bx^{(n)},q_n},$$
$$R_n(s) := \sum_{i=1}^{m(n)} \left( x_i^{(n)}  
1_{(\xi_i^{q_n} \leq \,s)} - \frac{(x_i^{(n)})^2}{\sigma_2(\bx^{(n)})}  s \right), \mbox{ and } Y_n(s):= Z_n(s) - R_n(s), \ s\geq 0.$$
It is implicit in the notation that $\xi_i^{q_n}:=\xi_i^{(n)}/q_n$, where $\xi_i^{(n)}$ has exponential (rate $x_i^{(n)}$) distribution, and where $(\xi_i^{(n)})_i$ are independent over $i$, for each $n$.

Define $\Zb_n,\Rb_n,\Yb_n$ to be respectively $Z_n,R_n,Y_n$ multiplied by $\frac1{\sigma_2(\bx^{(n)})}$, so that $\Zb_n\equiv \Yb_n+ \Rb_n$. 
It should not be surprising that both the shift in the \MC\ time and  the spatial scaling applied to the walks are the same as in \cite{EBMC}.
It is clear that, for each $n$, $R_n$ and $Y_n$ are independent (the former depends only of the first $m(n)$ $\xi^{(n)}$s, and the latter only on all the other $\xi^{(n)}$s).

Recall the definitions (\ref{deftilW}--\ref{defBtc}). The following result is a direct analogue of  \cite{EBMC}, Proposition 9.
\begin{Proposition}
\label{Panalognew}
If $(\kappa,0,\bc) \in {\cal I}$, and provided (\ref{hyp1}--\ref{hyp4}) are satisfied 
as $n \to \infty$, then 
$$(\Yb_n,\Rb_n) \cd ({\widetilde W}^{\kappa, t},V^\bc), \mbox{ as } n\to \infty,$$ 
where ${\widetilde
  W}^{\kappa, t}$ and $V^\bc$ are independent, and therefore $\Zb_n \cd W^{\kappa, t, \bc}.$
\end{Proposition}

The rest of this section is devoted to the proof of the above proposition, and some of its consequences.
As already mentioned, the argument is a simplification of that given in Section 2.4 of \cite{EBMC}, for the main reason that the current $Y_n$ has a simpler explicit form.
From now on assume that $t\in \mathbb{R}$ is given in Proposition \ref{Panalognew}.

Note that the independence of $Y_n$ and $R_n$ clearly implies that of $\Yb_n$ and $\Rb_n$. So provided that each of the sequences converges in law, the joint convergence in law to the product limit law is a trivial consequence. 
Furthermore, the convergence of $\Rb_n$ can be verified in a standard way (for each $k$, the $k$th largest jump of $\Rb_n$ converges to the $k$th largest jump of $V^\bc$, and the second-moment MG estimates are used to bound the tails), as was already done in \cite{EBMC}. 
Indeed, a careful look at the pages 17--19 in \cite{EBMC} suffices to see that the $\xi$s corresponding to the leading blocks were ``artificially reintegrated''  into consideration (see the last paragraph of \cite{EBMC}, page 17) via the random family $(\tilde{\xi}_i)_i$ which has exactly the same law as the present $(\xi_i^{q_n})_i$.
Therefore, the sequence of  processes appearing on the left-hand side of \cite{EBMC}, display (36) has the same law as the sequence $(\Rb_n)_n$, and we can simply restate that auxiliary result from \cite{EBMC} as 
 \begin{equation}
\label{Eimp1}
\Rb_n\cd V^\bc(s), \mbox{ as } n\to \infty.
\end{equation}

It remains to study the convergence of $(\Yb_n)_n$.  This sequence of processes differs from the equally named sequence in \cite{EBMC}.
Write $\sigma_r^n$ for $\sigma_r(\bx^{(n)})$, $r=1,2,3$ in the sequel.
An important observation is that 
\begin{equation}
\label{defYbn}
\Yb_n(s):= \sum_{i=m(n)+1}^{\len(\bx^{(n)})} \frac{x_i^{(n)}}{\sigma_2^n} 1_{(\xi^{q_n}_i\, \leq \ s)} -\frac{s}{\sigma_2^n} + 
\sum_{i=1}^{m(n)} \frac{(x_i^{(n)})^2}{(\sigma_2^n)^2}  s.
\end{equation}
The infinitesimal drift and variance calculations are now straightforward.
Let $\FF_s^n:=\sigma\{\Yb_n(u):u\leq s\}$, so that $\FF^n:=(\FF _s^n)_{s\geq 0}$ is the filtration generated by $\Yb_n$.
Using (\ref{defYbn}) one can easily derive that $E(d\Yb_n(s)|\,\FF_s^n) =$
\begin{equation}
 \left(\sum_{i=m(n)+1}^{\len(\bx^{(n)})} \frac{(x_i^{(n)})^2}{\sigma_2^n} \left(\frac1{\sigma_2^n} +t\right) 1_{(\xi^{q_n}_i > \ s)} - \frac{1}{\sigma_2^n} + \sum_{i=1}^{m(n)} \frac{(x_i^{(n)})^2}{(\sigma_2^n)^2} \right)ds=\nonumber
\end{equation}
\begin{eqnarray}
 \left(\sum_{i=m(n)+1}^{\len(\bx^{(n)})} \frac{(x_i^{(n)})^2}{\sigma_2^n} \left(\frac1{\sigma_2^n} +t\right) (1- 1_{(\xi^{q_n}_i \leq \ s)}) - \frac{1}{\sigma_2^n} + \sum_{i=1}^{m(n)} \frac{(x_i^{(n)})^2}{(\sigma_2^n)^2} \right)ds=\nonumber\\
\left(\sum_{i=m(n)+1}^{\len(\bx^{(n)})}\left[ \frac{(x_i^{(n)})^2}{\sigma_2^n}\,  (t- t1_{(\xi^{q_n}_i \leq \ s)}) - \frac{(x_i^{(n)})^2}{(\sigma_2^n)^2} 1_{(\xi^{q_n}_i \leq \ s)}) \right]\right)ds.
\label{Edrift1}
\end{eqnarray}
Let us estimate the remaining three sums separately.
In order to do so, note that (\ref{mChoo1new}--\ref{mChoo2new}) clearly imply 
$$
\sum_{i=m(n)+1}^{\len(\bx^{(n)})} \frac{(x_i^{(n)})^2}{\sigma_2^n} \to 1, \mbox{ as } n\to \infty.
$$
For the other two terms, use the approximation by the average. More precisely, the mean of the (absolute value of the) second sum can be bounded above by 
$$
t \sum_{i=m(n)+1}^{\len(\bx^{(n)})} \frac{(x_i^{(n)})^3}{\sigma_2^n} \cdot s \,q_n,
$$
and due to (\ref{hyp1},\ref{hyp4}), the last quantity becomes negligible as $n\to \infty$.
Finally note that it is elementary that if $a_{i,n}\geq 0$, then
$$\left(a_{i,n}1_{(\xi^{q_n}_i \leq s)} - a_{i,n} x_i^{(n)}  q_n s,\, s\geq 0\right)$$
is a supermartingale with Doob-Meyer decomposition
$M_{i,n}(s) - a_{i,n}(s-\xi_i^{q_n})^+$,
where $M_{i,n}$ is a martingale, and 
$\langle M_{i,n} \rangle (s) = a_{i,n}^2 x_i^{(n)}  q_n \min(s,\xi_i^{q_n})$
is its quadratic variation.
Now let $a_{i,n}:= \sfrac{(x_i^{(n)})^2}{(\sigma_2^n)^2}$, and note that due to (\ref{hyp1},\ref{hyp4},\ref{mChoo1new})
\begin{equation}
\label{Easykappa}
\sum_{i=m(n)+1}^{\len(\bx^{(n)})} a_{i,n} x_i^{(n)} q_n =\sum_{i=m(n)+1}^{\len(\bx^{(n)})} \frac{(x_i^{(n)})^3}{(\sigma_2^n)^2} q_n \to \kappa, \mbox{ as } n\to \infty,
\end{equation}
while 
$$
\sum_{i=m(n)+1}^{\len(\bx^{(n)})}  (a_{i,n})^2 = \sum_{i=m(n)+1}^{\len(\bx^{(n)})}  \frac{(x_i^{(n)})^4}{(\sigma_2^n)^4}
\leq \frac{x_{m(n)}^{(n)}\sigma_3^n}{(\sigma_2^n)^4},
$$
and 
$$
\sum_{i=m(n)+1}^{\len(\bx^{(n)})}  (a_{i,n})^2  x_i^{(n)} q_n = \sum_{i=m(n)+1}^{\len(\bx^{(n)})}  \frac{(x_i^{(n)})^5}{(\sigma_2^n)^4} (t+ 1/\sigma_2^n) 
=O\left( \frac{(x_{m(n)}^{(n)})^2\sigma_3^n}{(\sigma_2^n)^5}\right),
$$
both become negligible as $n\to \infty$ due to (\ref{hyp1}, \ref{mChoo1new}) and the fact that $\bc \in \cvt$.

Combining the just made observations with the usual $L^2$ (Doob) martingale estimates, one can 
conclude from (\ref{Edrift1}) that, for each fixed $s$,
\begin{equation}
\label{Edrift2}
E(d\Yb_n(s)|\,\FF_s^n) \cp (t - \kappa s)\,ds, \mbox{ as }n\to \infty.
\end{equation}
Recalling the representation (\ref{defYbn}), one can even more easily verify (via  (\ref{Easykappa}) given above) that, for each fixed $s$,
\begin{equation}
\label{Eivar}
E((d\Yb_n(s))^2|\,\FF_s^n)\cp \kappa \,ds, \mbox{ as }n\to \infty.
\end{equation} 
Since the largest jump of $\Yb_n$ is of size $x_{m(n)+1}^{(n)}/\sigma_2^n=o_n(1)$, 
the classical martingale central limit theorem (cf.~\cite{ethierkurtz}) implies that 
\[
\Yb_n \cd {\widetilde W}^{\kappa, t}, \mbox{ as }n\to \infty,
\]  
and, as already argued, this concludes the proof of the proposition.

{\em Remark.} By the Cauchy-Schwarz inequality,
$(\sigma_2^n)^2 \leq \sigma_1^n \sigma_3^n$,
and so (\ref{hyp1},\ref{hyp4}) imply that $\sigma_1^n\to \infty$ as $n\to \infty$.
While this fact was needed in the proof of the analogous \cite{EBMC}, Proposition 9, here it could slip by unnoticed. 
If $\kappa>0$, it is easy to see that the limit $W^{\kappa, t,\bc}$ of $\Zb_n$ has (countably) infinitely many excursions above past minima. If $\kappa=0$ and $\bc \in \cvt \setminus \cvd$, the same was proved in \cite{EBMC}, Proposition 14.
\hfill\endrem

Using the Skorokhod representation theorem, we may assume that the convergence stated in Proposition \ref{Panalognew} holds in the almost sure sense. 
To state the next result (essential for the conclusions to be made in Section \ref{S:C}), redefine $q_n(t):= t+ \frac{1}{\sigma_2^n}$, for $t\in \mathbb{R}$.
Then let $Z_n^t:=Z^{\bx^{(n)},q_n(t)}$ and $\bar{Z}_n^t:= Z_n^t/\sigma_2^n$.
The (almost sure version of) Proposition \ref{Panalognew} says that there exists a Brownian motion $W$ and an independent  jump process $V^\bc$, such that 
$$
\bar{Z}_n^t \to W^{\kappa,t,\bc}, \mbox{ almost surely, as } n\to \infty,
$$
where the convergence of paths is considered in the Skorokhod $J_1$ topology.
Let 
$$
A_t:=\{\omega: \bar{Z}_n^t(\cdot)(\omega) \to W^{\kappa,t,\bc}(\cdot)(\omega) \mbox{ in the Skorokhod $J_1$ topology}\}.
$$
\begin{Lemma}
\label{Lsimultcon}
On the event $A_t$, for any $z\in \mathbb{R}$
$$
(\bar{Z}_n^z(s),\,s\geq 0) \to (W^{\kappa,t,\bc}(s) + (z-t)s,\,s\geq 0) \equiv W^{\kappa,z,\bc}, \mbox{ as } n\to \infty,
$$
in the Skorokhod $J_1$ topology.
\end{Lemma}
{\em Proof.} Recall the explicit form (\ref{defZbxq}) of $Z^{\,\cdot,\cdot}$
Observe the following identity:
\[
Z_n^z\left( s \cdot \frac{q_n(t)}{q_n(z)}\right) = Z_n^t(s) + s \left( 1- \frac{q_n(t)}{q_n(z)}\right), \ \forall s\geq 0.
\]
Since clearly $\lim_{n\to \infty} \frac{q_n(t)}{q_n(z)}= 1$, and moreover since
\[
\frac{1}{\sigma_2^n} \left( 1- \frac{q_n(t)}{q_n(z)}\right) = z - t + O_{z,t}(\sigma_2^n)\,,
\]
the convergence stated in the lemma follows omega-by-omega on $A_t$. 

\section{Conclusions}
\label{S:C}
Here is an immediate consequence of Lemma \ref{Ltosimplify} and \cite{EBMC}, Proposition 7 and 9, 
which was announced in Remark (d) of Section \ref{S:SBFW}.
\begin{Corollary}
\label{C:immediate}
For each fixed $t$, under the hypotheses of Proposition \ref{Panalognew}, the sequence
$(\bX^{(n)}(q_n(t)))_n$ converges in law (with respect to $\cvd$-metric) to the sequence of ordered excursions of $B^{\kappa,t,\bc}$ (away from $0$) .
\end{Corollary}

But in fact more is true in view of Lemma \ref{Lsimultcon}.
From now on we take the families $(W^{\kappa,t,\bc},\,t\in \mathbb{R})$ and  $(B^{\kappa,t,\bc},\,t\in \mathbb{R})$
as jointly defined, on a common probability space, via a given pair $(W,V^\bc)$, where $W$ is Brownian motion and $V^\bc$ is an independent jump process from (\ref{defVc}).

Let us denote by $A$ the event $A_t$ of full probability from Lemma \ref{Lsimultcon}. 
For each $t\in \mathbb{R}$, define $\Xi^{(n)}(t)$  to be the point process on $[0,\infty)\times (0,\infty)$ such that $(x,y)$ is in $\Xi^{(n)}(t)$ if and only if there is an excursion above past minima of $\bar{Z}_n^t$ (see (\ref{Eexcur}) and Figure 2), starting from $x$ and ending at $x+y$.
Similarly, let $\Xi^{(\infty)}(t)$  be the point process on $[0,\infty)\times (0,\infty)$ such that $(x,y)$ is in $\Xi^{(\infty)}(t)$ if and only if there is an excursion away from $0$ of $B^{\kappa,t,\bc}$, starting from $x$ and ending at $x+y$.
One can then apply deterministic result stated as \cite{aldRGMC}, Lemma 7 to conclude the following: on the event $A$ of full probability,
 for each $t\in \mathbb{R}$, one has 
\begin{equation}
\label{EconvXis}
\lim_{n\to \infty}\Xi^{(n)}(t) = \Xi^{(\infty)}(t),
\end{equation}
in the sense of vague convergence of counting measures on
$[0,\infty) \times (0,\infty)$ (see e.g.~\cite{kal83}).
As in \cite{aldRGMC,EBMC}, write $\pi$ for the ``project onto the $y$-axis'' defined on $\mathbb{R}^2$, and ``${\rm ord}$'' for the ``decreasing
ordering" map defined on infinite-length vectors, respectively.
For a fixed (think large) $K<\infty$, 
define in addition
$\pi_K$ to be the ``project the strip $[0,K]\times (0,\infty)$ onto the $y$-axis'' analogue of $\pi$. 
Then 
one can recognize ${\rm ord}(\pi(\Xi^{(n)}(t)))$ as $\bX^{(n)}(q_n(t))$, and 
$\bX^{(\infty)}(t):={\rm ord}(\pi(\Xi^{(\infty)}(t)))$ as the infinite vector of excursion lengths of $B^{\kappa,t,\bc}$.
Similarly $\pi_K(\Xi^{(\infty)}(t))$ (resp.~$\pi_K(\Xi^{(n)}(t))$) is the collection of all the excursions of $B^{\kappa,t,\bc}$ (resp.~$\bar{Z}_n^t$), which start before time $K$. 

We already know that the law of $\bX^{(\infty)}(t)$ is that of  the marginal 
of $\mu(\kappa,0,\bc)$ at time $t$ (or equivalently, the marginal of $\mu(\kappa,t,\bc)$ at time $0$).
The vague convergence (\ref{EconvXis}) now easily implies that 
there exists a (random) order of $\pi_K(\Xi^{(n)}(t))$, here temporarily denoted by ${\rm ord}_{\Xi^{(n)},\Xi^{(\infty)}}$, since it is 
induced by the similarity of the $\Xi$s, such that 
\begin{equation}
\label{Edordnord}
\|{\rm ord}_{\Xi^{(n)},\Xi^{(\infty)}} (\pi_K(\Xi^{(n)}(t)))- {\rm ord}(\pi_K(\Xi^{(\infty)}(t))) \|_2 \to 0,  \mbox{ on the event $A$}.
\end{equation}
In words, if considering only the starts before time $K$, it is possible to order the excursions of $\bar{Z}_n^t$ so that 
 the corresponding infinite vector (obtained by appending an infinite sequence of $0$s to the elements of $\pi_K(\Xi^{(n)}(t))$) 
 matches the infinite vector $ {\rm ord}(\pi_K(\Xi^{(\infty)}(t)))$ in $l_1$-norm 
up to an $o_n(1)$ error term. Moreover, convergence in $l_1$-norm implies convergence in $l_2$-norm.

Take $z>t$, another real number, and let $\eps>0$ be fixed but arbitrarily small.
From now on $u$ will denote either $t$ or $z$.
In order to upgrade (\ref{Edordnord}) to the convergence of $\bX^{(n)}(q_n(u))$ to $\bX^{(\infty)}(u)$ with respect to distance $d(\cdot,\cdot)$, one can make the following observations.
For $\bx \in \cvd$, let $f_m(\bx) := (x_1,\ldots,x_m,0,0,\ldots)$ be the ``projection'' onto the first $m$ components. 
Take some arbitrarily large integer $k$, and choose
$m_k\in \mathbb{N}$ such that 
\begin{equation}
\label{Econtr1a}
P(d(\bX^{(\infty)}(u),f_{m_k}(\bX^{(\infty)}(u))) >\eps)<\frac{1}{2^k}.
\end{equation}
Since $\bX^{(n)}(q_n(t))\cd \bX^{(\infty)}(t)$ with respect to $d$, for this $k$ and possibly larger but still finite $m_k$ we can have in addition
\begin{equation}
\label{Econtr1b}
\limsup_n P(d(\bX^{(n)}(q_n(u)),f_{m_k}(\bX^{(n)}(q_n(u)))) >\eps)< \frac{1}{2^k}.
\end{equation}
In words, with an overwhelming probability, all the (random) infinite vectors under consideration are well-approximated (in the $l_2$-norm) by their first $m_k$ components.

Since $\pi(\Xi^{(\infty)}(\cdot))$ is $l_2$-valued, and since its elements are listed in size-biased order, one can easily deduce that 
for the above $m_k$, there exists some large time $K_k:=K(m_k)<\infty$, such that 
\begin{equation}
\label{Econtr2a}
P(f_{m_k}({\rm ord}(\pi_K(\Xi^{(\infty)}(u)))) \neq f_{m_k}(\bX^{(\infty)}(u)))<\frac{1}{2^k}.
\end{equation}
In words, $K$ is sufficiently large so that with high probability the largest $m_k$ elements of $\pi(\Xi^{(\infty)}(u))$,  all correspond to excursions that started before time $K$.
Again due to $\bX^{(n)}(q_n(u))\cd \bX^{(\infty)}(u)$, the analogous 
\begin{equation}
\label{Econtr2b}
\limsup_n P(f_{m_k}({\rm ord}(\pi_K(\Xi^{(n)}(q_n(u))))) \neq f_{m_k}(\bX^{(n)}(q_n(u))))<\frac{1}{2^k}
\end{equation}
is implied for some (possibly larger but) finite $K=K(m_k)$.

Apply the triangle inequality to bound $d(\bX^{(n)}(q_n(u)), \bX^{(\infty)}(u))$  by the sum of three terms: $d(\bX^{(n)}(q_n(u)), f_{m_k}(\bX^{(n)}(q_n(u))))$, 
$d(f_{m_k}(\bX^{(n)}(q_n(u))), f_{m_k}(\bX^{(\infty)}(u)))$, and $d(f_{m_k}(\bX^{(\infty)}(u)), \bX^{(\infty)}(u))$.
The initial and the final term are controlled by (\ref{Econtr1a}--\ref{Econtr1b}), while the middle term is controlled by (\ref{Econtr2a}--\ref{Econtr2b}) and (\ref{Edordnord}), where one makes use of the elementary inequality: for $\bx,\by \in l^2$,
$$
	d({\rm ord}(\bx), {\rm ord}(\by)) \leq \sum_i (x_i - y_{b(i)})^2,
$$
regardless of the choice of bijection $b:\mathbb{N}\to \mathbb{N}$. 

{\em Remark.} It is clear (for example from (\ref{Edrift2}--\ref{Eivar}), think about redefining $q_n(t)$ as 
$\sfrac{1}{\sigma_2(\bx^{(n)})} + t-\tau$), that the parameter $\tau$ corresponds to the time-shift of the eternal \MC, and so the above conclusions automatically extend to the setting where $\tau\neq 0$. 
\hfill\endrem

An immediate conclusion is  
\begin{Lemma}
\label{LjointconX}
If 
$
((\Zb ^{x^{(n)},q_n(t)}(s),\,s \geq - 1/\sigma_2^n),\,t \in \mathbb{R}) \longrightarrow ((W^{\kappa,t-\tau,\bc}(s),\,s\geq 0)\,, t \in \mathbb{R}),
\mbox{ as } n\to \infty,
$
in the sense of Lemma \ref{Lsimultcon},  and
 if  $\bX^{(\infty)}(t)\in\cvd$ is the vector of ordered excursion lengths of $B^{\kappa,t-\tau,\bc}$, then\\ 
(i) for any $t\in \mathbb{R}$ 
\[
d(\bX^{(n)}(q_n(t)), \bX^{(\infty)}(t)) \cp 0, \mbox{  as }n\to \infty,
\]
(ii) for any finite sequence of times $t_1<t_2,\ldots<t_m$, one can find a subsequence $(n_j)_j$ such that almost surely
\[
d(\bX^{(n_j)}(q_{n_j}(t_k)),\bX^{(\infty)}(t_k))\to 0, \mbox{ for all }  k=1,\ldots,m, \mbox{  as }j\to \infty.
\] 
\end{Lemma}


Recalling that  $(\bX^{(n)}(q_n(s)),\,s \geq - 1/\sigma_2^n)$ has the law of the \MC\ (see Proposition \ref{Pcoro}),
and applying the Feller property together with Lemma \ref{LjointconX}(i), where one should identify $\bX^{(\infty)}$ with $\bX$,
will  complete the proof of the claim about the distribution of $\bX$ in Theorem \ref{Tmain}.
It is easy to see (using arguments analogous to those given above) that the realization $\bX^{(\infty)}\equiv\bX$ of each eternal version from Theorem \ref{Tmain} is a c\`adl\`ag (rcll) process on an event of full probability.

{\em Remark.} Recall the {\rm COL} operation of \cite{EBMC}, Section 5. In particular, each $c_i>0$ is interpreted as the rate of Poisson coloring (per unit mass) by marks of the $i$th ``color'', applied to the standard Aldous' \MC\ $\bX^*$. Once all the color marks are deposited, any two blocks
of  $\bX^*$ that share at least one mark of the same color are instantaneously and simultaneously merged together. The jump in $W^{\kappa, \cdot -\tau, \bc}$ at time $\xi_i$ of size $c_i$ corresponds precisely to the effect of coloring by the $i$th color.
Moreover, one could argue that if $\widetilde{W}^{\kappa,\cdot -\tau,{\bf 0}}$ and $\widetilde{W}^{\kappa,\cdot-\tau,\bc}$ are given  in
(\ref{defWtc}) using the same Brownian motion $W$, then the excursions of the corresponding $(B^{\kappa,t-\tau+\|\bc\|_2,\bc},\,t\in \mathbb{R})$
(suppose for simplicity that $\bc\in \cvd$) away from $0$ are almost surely the result of the above ${\rm COL}$ operation executed on the excursions of $(B^{\kappa,t-\tau,{\bf 0}},\,t\in \mathbb{R})$.
The fact that, as time increases, each color ``spreads'' in this coupling  (almost surely) 
 only over the ``neighboring'' blocks may again seem counterintuitive.
The point is that ${\rm COL}$ commutes with the \MC\ dynamics, and that therefore it can be pushed to $-\infty$. 
The infinitesimally small dust particles of $\bX^*(-\infty)$ are mutually interchangeable. 
The $i$th color at time $-\infty$ is represented as an additional dust particle, of mass much superior to
standard dust, but still negligible (a formal statement of this is (\ref{Thyp3}) or (\ref{hyp3})).
One can naturally couple the representation of the \MC\ using simultaneous breadth-first walks (or Uribe's diagram)
started only from standard dust as $t\to -\infty$, with the same representation of the \MC\ started from the union of two types of dust as $t\to-\infty$. Proposition \ref{Panalognew} and Lemma \ref{Lsimultcon} do this formally (their predecessor is \cite{EBMC}, Proposition 41).
In this coupling every color gradually spreads only to neighboring blocks of those already marked by it.

\section{A list of questions with some remarks}
\label{S:Q}
{\bf 1.} Fix $\bc, \kappa$ and $\tau$ as in Theorem \ref{Tmain}.
The convergence of Lemma \ref{LjointconX} implies that of $(\bX^{(n)}(q_n(t)),\, t\in \mathbb{R})$ to 
$(\bX^{(\infty)}(t),\, t\in \mathbb{R})$ in the sense of f.d.d.~as $n\to \infty$. 
It seems highly plausible that this could be extended to the convergence in law with respect to the Skorokhod $J_1$ topology on $\cvd$-valued path processes, using the well-known Aldous' criterion for tightness, and bounds 
on the increase of the second-moment $\sigma_2(\bX^{(n)}(q_n(\cdot)))$ over small intervals of time, as in \cite{aldRGMC,EBMC}.\\
A more challenging/interesting project would be to try to strengthen the convergence of Lemma \ref{LjointconX}(ii) to the full convergence in the almost sure sense. Having that, could one prove that on some set of full probability we have
$$(\bX^{(n)}(q_n(t)),\, t\in \mathbb{R})\to  (\bX^{(\infty)}(t),\, t\in \mathbb{R})$$
with respect to the Skorokhod $J_1$ distance on $D((-\infty,\infty),\cvd)$?\\
{\bf 2.} The author tried a ``lazy'' approach of comparing directly, for each fixed large $n$, the current process $\Zb_n$ to the analogous $\Zb_n$ from \cite{EBMC}, in order to be able to conclude Proposition \ref{Panalognew} without any additional effort. The absence of the load-free intervals for $\Zb_n$ from \cite{EBMC} makes the direct ($J_1$ distance) comparison of the paths of two pre-limit processes difficult, unless one uses the extra information that they both converge to the same limit in the Skorokhod $J_1$ sense. Perhaps there is a trick, which could make the direct comparison possible, so that the conclusion about the same distributional limit could be made in a handwaving manner?\\
 {\bf 3.} Is there a direct connection between Prim's algorithm of \cite{bromar15} and the simultaneous breadth-first walks of Section \ref{S:SBFW}? The work in progress \cite{uribe_wip}, joint with the observations of Section \ref{S:BFWUD}, could provide a two-step connection between the two representations. \\
{\bf 4.} 
Theorem \ref{Tmain} enables one to construct a coupling ${\mathcal X}$ of the whole family $(\bX^{\kappa,\tau,\bc}: (\kappa,\tau,\bc)\in {\cal I})$ on one and the same probability space, using a single Brownian motion $W$ (in order to define $\widetilde{W}^{\cdot,\cdot}$ as in (\ref{deftilW})) and an independent countable family of i.i.d. exponential (rate $1$) random variables (in order to define $V^{\cdot}$ as in (\ref{defVc})). Here it is natural to follow the convention that if $c_i=0$ for some $i$, then the corresponding exponential divided by $c_i$ is almost surely infinite, and can therefore be omitted from the sum in $V^{\bc}$. In particular $V^{{\bf 0}}$ is the zero process. The family of deterministic constant coalescents $\hat{\mu}(y)$, $y\geq 0$ can clearly be included in the above coupling ${\mathcal X}$. It is not hard to see that on an event of full probability all the elements of ${\mathcal X}$ have c\`adl\`ag paths.
Is it possible to rephrase the questions asked in the final remarks (Section 6) of \cite{EBMC} in terms of the almost sure distance on $D((-\infty,\infty),\cvd)$ between the elements of $\mathcal X$.
In particular, does the sequence of laws $(\mu(\kappa_m,\tau_m,\bc_m))_m$ converge if and only if $(\bX^{\kappa_m,\tau_m,\bc_m})_m$
converges in $D((-\infty,\infty),\cvd)$, and if and only if $(B^{\kappa_m,\tau_m,\bc_m})_m$ converges in $D((-\infty,\infty),[0,\infty))$?
Is ${\mathcal X}$ closed in $D((-\infty,\infty),\cvd)$?
Which sequences $(\bX^{\kappa_m,\tau_m,\bc_m})_m$ in ${\mathcal X}$ converge to the zero process $V^{{\bf 0}}$?
\\ 
{\bf 5.}
The fragmentation/coalescence duality announced in the title of \cite{armen05} has a chance of being explored again in the more general context of Theorem \ref{Tmain}. 
What are the fragmentation rates of the process dual to each $\bX$ of Theorem \ref{Tmain}?
Could the dual be used to derive new properties of the \MC\ entrance boundary? \\
For example. take $\bX^{\kappa,\tau,\bc}(t)$, where $\bX^{\kappa,\tau,\bc}$ has distribution  $\mu(\kappa,\tau,\bc)$.
Equivalently, let  $\bX^{\kappa,\tau,\bc}(t)$ be the vector of ordered excursion (away from $0$) lengths of $B^{\kappa,\tau,\bc}$.
We know that
\begin{equation}
 {\cal L}(\bX^{\kappa,\tau,\bc}(t)) = \int_{{\cal I}} 
{\cal L}(\bX^{\cdot,\cdot,\cdot}(t)) d\nu(\cdot,\cdot,\cdot) , \ 
\label{mixture}
\end{equation}
where $\nu$ is the Dirac mass at $(\kappa,\tau,\bc)$.
Can (\ref{mixture}) hold for another (non-trivial) measure on ${\cal I}$? Equivalently, are already the marginal laws $ {\cal L}(\bX^{\cdot,\cdot,\cdot}(t))$  mutually singular, or could $\bX^{\kappa_1,\tau_1,\bc_1}(t)$ and $\bX^{\kappa_2,\tau_2,\bc_2}(t)$ be absolutely continuous?
In the latter case, what is the corresponding Radon-Nikod\'ym derivative?\\
{\bf 6.} The asymptotic behavior of the diagram $L(\xi^{(n)},\bx^{(n)})$  and its related ``genealogical tree''-analog $T$ (see also Figures 3 and 4) of Section \ref{S:UD} could be related to the previous question. Under hypotheses (\ref{hyp1}--\ref{hyp4}) on $\bx^{(n)}$, the sequence of corresponding diagrams/trees should formally converge to a tree-like object that encodes the corresponding limiting extreme eternal version of the \MC. Can this be formalized? If so, does the limiting tree have anything in common with the ICRT of \cite{bhahof3}? 
Addario-Berry et al.~\cite{abbgm2} work on minimal spanning tree convergence to the above tilted ICRT and show that the limiting object has fractal dimension $3$. Does this help in answering the above question?\\
{\bf 7.} Could the definition of Uribe's diagram be extended and usefully exploited in the context of more general merging kernels? 
\\
{\bf 8.} Does there exist a non-standard augmented \MC, analogous to the standard augmented \MC\ of Bhamidi et al? Is it the scaling-limit 
for the coupled random graph and the corresponding excess-edge data count process
under hypotheses (\ref{hyp1}-\ref{hyp4})?
Addario-Berry et al.~\cite{abbgm1} 
construct of a ``metric multiplicative coalescent'' in the Gromov-Hausdorff-Prokhorov topology, which may provide a good framework for making progress in this direction.
Could \cite{edgsur} somehow complement that?
\\
{\bf 9.} Suppose that $\bc\neq {\bf 0}$ and let the family $(W^{\kappa,t-\tau,\bc},B^{\kappa,t-\tau,\bc})$, $t\in \mathbb{R}$ be coupled as in Theorem \ref{Tmain}. The variance scale $\kappa$ could be positive or $0$. 
Using the \MC\ representation of Theorem \ref{Tmain} from the viewpoint of of  the ``coloring construction''  (see 
\cite{EBMC}, Section 5 and the final remark of Section \ref{S:C}) it is immediate that, almost surely, 
simultaneously for all $t\in \mathbb{R}$, no (positive) jump of $V^{\bc}$ arrives at the very beginning of any excursion of $B^{\kappa,t-\tau,\bc}$ away from $0$. 
In other words, no excursion of $B^{\kappa,\cdot-\tau,\bc}$ away from $0$ ever starts with a jump.
Is this obvious from the stochastic calculus point of view, and why?

\end{document}